\documentclass[a4paper,12pt]{article}
\usepackage{amsmath}
\usepackage{amscd}
\usepackage[xdvi]{graphicx}
\usepackage[dvips]{color}
\usepackage{times}
\usepackage{txfonts}

\makeatletter

\@addtoreset{equation}{section}

\newcounter{def}[section]
\renewcommand{\thedef}{\stepcounter{def}\thesection.\@arabic\c@def }

\begin{document}
\setlength{\baselineskip}{24pt}
\begin{center}
\textbf{\LARGE{A spectral theory of linear operators on rigged Hilbert spaces under analyticity conditions II:
applications to Schr\"{o}dinger operators}}
\end{center}

\setlength{\baselineskip}{14pt}

\begin{center}
Institute of Mathematics for Industry, Kyushu University, Fukuoka,
819-0395, Japan

\large{Hayato CHIBA} \footnote{E mail address : chiba@imi.kyushu-u.ac.jp}
\end{center}
\begin{center}
May 20, 2015
\end{center}

\begin{center}
\textbf{Abstract}
\end{center}
A spectral theory of linear operators on a rigged Hilbert space is applied to Schr\"{o}dinger operators with
exponentially decaying potentials and dilation analytic potentials.
The theory of rigged Hilbert spaces provides a unified approach to resonances (generalized eigenvalues)
for both classes of potentials without using any spectral deformation techniques.
Generalized eigenvalues for one dimensional Schr\"{o}dinger operators (ordinary differential operators) are investigated in detail.
A certain holomorphic function $\mathbb{D}(\lambda )$ is constructed so that $\mathbb{D}(\lambda ) = 0$
if and only if $\lambda $ is a generalized eigenvalue.
It is proved that $\mathbb{D}(\lambda )$ is equivalent to the analytic continuation of the Evans function.
In particular, a new formulation of the Evans function and its analytic continuation is given.
\\[0.2cm]
\textbf{Keywords}: spectral theory; resonance pole; rigged Hilbert space; generalized spectrum; Schr\"{o}dinger operator



\section{Introduction}

A spectral theory of linear operators is one of the fundamental tools in functional analysis 
and well developed so far.
Spectra of linear operators provide us with much information about the operators
such as the asymptotic behavior of solutions of linear differential equations.
However, there are many phenomena that are not explained by spectra.
For example, transient behavior of solutions of differential equations is not described by spectra;
even if a linear operator $T$ does not have spectrum on the left half plane,
a solution of the linear evolution equation $dx/dt = Tx$ on an infinite dimensional space
can decay exponentially as $t$ increases for a finite time interval.
Now it is known that such transient behavior can be induced by resonance poles or generalized eigenvalues,
and it is often observed in infinite dimensional systems such as plasma physics \cite{Cra},
coupled oscillators \cite{Chi2,Str2} and Schr\"{o}dinger equations \cite{His, Reed}.

In the literature, resonance poles for Schr\"{o}dinger operators $-\Delta + V$ are defined in several ways.
When a wave operator and a scattering matrix can be defined, resonance poles may be defined as poles of 
an analytic continuation of a scattering matrix \cite{Reed}.
When a potential $V(x)$ decays exponentially, resonance poles can be defined with the aid of certain weighted 
Lebesgue spaces, which is essentially based on the theory of rigged Hilbert spaces \cite{Rau}.
When a potential has an analytic continuation to a sector around the real axis,
spectral deformation (complex distortion) techniques are often applied to define resonance poles,
see \cite{His} and references therein.

The spectral theory based on rigged Hilbert spaces (Gelfand triplets) was introduced by Gelfand et al.\cite{Gel1}
to give generalized eigenfunction expansions of selfadjoint operators.
Although they did not treat with resonance poles, the spectral theory of resonance poles (generalized spectrum)
of selfadjoint operators based on rigged Hilbert spaces are established by Chiba \cite{Chi} without 
using any spectral deformation techniques.

Let $\mathcal{H}$ be a Hilbert space, $X$ a topological vector space, which is densely and continuously 
embedded in $\mathcal{H}$, and $X'$ a dual space of $X$.
A Gelfand triplet (rigged Hilbert space) consists of three spaces $X \subset \mathcal{H} \subset X'$.
Let $T$ be a selfadjoint operator densely defined on $\mathcal{H}$.
The resolvent $(\lambda -T)^{-1}$ exists and is holomorphic on the lower half plane, while
it does not exist when $\lambda $ lies on the spectrum set $\sigma (T) \subset \mathbf{R}$.
However, for a ``good" function $\phi$, $(\lambda -T)^{-1}\phi$ may exist on $\sigma (T)$ in some sense and
it may have an analytic continuation from the lower half plane to the upper half plane
by crossing the continuous spectrum on the real axis.
The space $X$ consists of such good functions with a suitable topology.
Indeed, under certain analyticity conditions given in Sec.2, it is shown in \cite{Chi} that the resolvent has an analytic continuation
from the lower half plane to the upper half plane, which is called the generalized resolvent $\mathcal{R}_\lambda $ of $T$,
even when $T$ has the continuous spectrum on the real axis.
The generalized resolvent is a continuous operator from $X$ into $X'$, 
and it is defined on a nontrivial Riemann surface of $\lambda $.
The set of singularities of $\mathcal{R}_\lambda $ on the Riemann surface is called the generalized spectrum of $T$.
The generalized spectrum consists of a generalized point spectrum, a generalized continuous spectrum and 
a generalized residual spectrum set, which are defined in a similar manner to the usual spectral theory.
In particular, a point $\lambda $ of the generalized point spectrum is called a generalized eigenvalue.
If a generalized eigenvalue is not an eigenvalue of $T$ in the usual sense, it is called a resonance pole 
in the study of Schr\"{o}dinger operators.
A generalized eigenfunction, a generalized eigenspace and the multiplicity associated with a generalized
eigenvalue are also defined.
The generalized Riesz projection $\Pi$ is defined through a contour integral of $\mathcal{R}_\lambda $ as usual.
In \cite{Chi}, it is shown that they have the same properties as the usual theory.
For example, the range of the generalized Riesz projection $\Pi$ around an isolated generalized eigenvalue
coincides with its generalized eigenspace.
Although this property is well known in the usual spectral theory, our result is nontrivial because
$\mathcal{R}_\lambda $ and $\Pi$ are operators from $X$ into $X'$, so that
the resolvent equation and the property of the composition $\Pi \circ \Pi = \Pi$ do not hold.
If the operator $T$ satisfies a certain compactness condition, the Riesz-Schauder theory on a rigged Hilbert space
is applied to conclude that the generalized spectrum consists of a countable number of generalized eigenvalues
having finite multiplicities.
It is remarkable that even if the operator $T$ has the continuous spectrum (in the usual sense),
the generalized spectrum consists only of a countable number of generalized eigenvalues
when $T$ satisfies the compactness condition.

In much literature, resonance poles are defined by the spectral deformation techniques.
The formulation of resonance poles based on a rigged Hilbert space has the advantage that 
generalized eigenfunctions, generalized eigenspaces and the generalized Riesz projections 
associated with resonance poles are well defined and
they have the same properties as the usual spectral theory, although in the formulation based on the
spectral deformation technique, correct eigenfunctions associated with resonance poles 
of a given operator $T$ is not defined because $T$ itself is deformed by some transformation.
The defect of our approach based on a rigged Hilbert space is that a suitable topological vector space $X$
has to be defined, while in the formulation based on the spectral deformation technique,
a topology need not be introduced on $X$ because resonance poles are defined by using the deformed operator
on the Hilbert space $\mathcal{H}$, not $X$.
Once the generalized eigenfunctions and the generalized Riesz projections associated with resonance poles
are obtained, they can be applied to the dynamical systems theory.
The generalized Riesz projection for an isolated resonance pole on the left half plane (resp. on the imaginary axis)
gives a stable subspace (resp. a center subspace) in the generalized sense.
They are applicable to the stability and bifurcation theory \cite{Chi2} involving essential spectrum on the imaginary axis. 

In this paper, the spectral theory based on a rigged Hilbert space is applied to Schr\"{o}dinger operators
$T = -\Delta + V$ on $L^2(\mathbf{R}^m)$, where $\Delta$ is the Laplace operator and $V$ is the multiplication 
operator by a function $V(x)$.
Two classes of $V$ will be considered.
\\

(I) \textit{Exponentially decaying potentials.}
Suppose that $V$ satisfies $e^{2a|x|}V(x) \in L^2(\mathbf{R}^m)$ for some $a>0$.
Then, a suitable rigged Hilbert space for $T = -\Delta + V$ is given by
\begin{equation}
L^2(\mathbf{R}^m, e^{2a|x|}dx) \subset L^2(\mathbf{R}^m) \subset L^2(\mathbf{R}^m, e^{-2a|x|}dx).
\end{equation}
When $m$ is an odd integer, 
it is proved that the resolvent $(\lambda -T)^{-1}$ has a meromorphic continuation to the Riemann surface defined by
\begin{equation}
P(a) = \{ \lambda \, | \, -a < \mathrm{Im} (\sqrt{\lambda } ) < a \},
\end{equation}
(see Fig.\ref{fig2}) as an operator from $L^2(\mathbf{R}^m, e^{2a|x|}dx)$ into $L^2(\mathbf{R}^m, e^{-2a|x|}dx)$.
When $m$ is an even integer, $(\lambda -T)^{-1}$ has a meromorphic continuation to a similar region
on the logarithmic Riemann surface.
\\

(II) \textit{Dilation analytic potentials.}
Suppose that $V \in G(-\alpha ,\alpha )$ for some $0<\alpha < \pi/2$, where $G(-\alpha ,\alpha )$
is the van Winter space consisting of holomorphic functions on the sector $\{ z\, | \, -\alpha < \mathrm{arg} (z) < \alpha \}$,
see Sec.3.2 for the precise definition.
Then, a suitable rigged Hilbert space for $T = -\Delta + V$ is given by
\begin{equation}
G(-\alpha ,\alpha ) \subset L^2(\mathbf{R}^m) \subset G(-\alpha ,\alpha )',
\end{equation}
when $m=1,2,3$.
In this case, 
it is proved that the resolvent $(\lambda -T)^{-1}$ has a meromorphic continuation to the Riemann surface defined by
\begin{equation}
\hat{\Omega } = \{ \lambda \, | \, -2\pi-2\alpha < \mathrm{arg} (\lambda ) < 2\alpha \},
\end{equation}
as an operator from $G(-\alpha , \alpha )$ into $G(-\alpha ,\alpha )'$.
\\

For both cases, we will show that $T = -\Delta + V$ satisfies all assumptions for our spectral theory given in Sec.2,
so that the generalized resolvent, spectrum and Riesz projection are well defined.
In particular, the compactness condition is fulfilled, which proves that the generalized spectrum
on the Riemann surface consists of a countable number of generalized eigenvalues.
This result may be well known for experts, however, we will give a unified and systematic approach for 
both classes of potentials.
In the literature, exponentially decaying potentials are investigated with the aid of the weighted Lebesgue
space $L^2(\mathbf{R}^m, e^{2a|x|}dx)$ \cite{Rau} as in the present paper, while for dilation analytic potentials,
the spectral deformations are mainly used \cite{His}.
In the present paper, resonance poles are formulated by means of rigged Hilbert spaces
without using any spectral deformations for both potentials.
Once resonance poles are formulated in such a unified approach, a theory of generalized spectrum developed in \cite{Chi}
is immediately applicable.
The formulation of resonance poles based on rigged Hilbert spaces plays a crucial role when applying it to
the dynamical systems theory because eigenspaces and Riesz projections to them are well defined 
as well as resonance poles.

In Sec.4, our theory is applied to one dimensional Sch\"{o}dinger operators
\begin{equation}
- \frac{d^2}{dx^2} + V(x),\quad x\in \mathbf{R}.
\label{1-3}
\end{equation}
For ordinary differential operators of this form, the Evans function $\mathbb{E}(\lambda )$ is often used to detect
the location of eigenvalues.
The Evans function for an exponentially decaying potential is defined as follows:
let $\mu_+$ and $\mu_-$ be solutions of the differential equation
\begin{equation}
\left( \frac{d^2}{dx^2} + \lambda - V(x) \right) \mu = 0
\label{1-4}
\end{equation}
satisfying the boundary conditions
\begin{equation}
\left(
\begin{array}{@{\,}c@{\,}}
\mu_+(x, \lambda ) \\
\mu_+'(x, \lambda )
\end{array}
\right) e^{\sqrt{-\lambda }x} \to \left(
\begin{array}{@{\,}c@{\,}}
1 \\
-\sqrt{-\lambda }
\end{array}
\right), \quad (x\to \infty), 
\label{1-1}
\end{equation}
and 
\begin{equation}
\left(
\begin{array}{@{\,}c@{\,}}
\mu_-(x, \lambda ) \\
\mu_-'(x, \lambda )
\end{array}
\right) e^{-\sqrt{-\lambda }x} \to \left(
\begin{array}{@{\,}c@{\,}}
1 \\
\sqrt{-\lambda }
\end{array}
\right), \quad (x\to -\infty), 
\label{1-2}
\end{equation}
respectively, where $\mathrm{Re} (\sqrt{-\lambda }) > 0$.
Then, $\mathbb{E}(\lambda )$ is defined to be the Wronskian
\begin{equation}
\mathbb{E}(\lambda ) = \mu_+(x, \lambda )\mu_-'(x, \lambda ) - \mu_+'(x, \lambda )\mu_-(x, \lambda ).
\end{equation}
It is known that $\mathbb{E}(\lambda )$ is holomorphic on $\{ \lambda \, | \, -2\pi < \mathrm{arg} (\lambda ) < 0\}$
(that is, outside the essential spectrum of $-d^2/dx^2$) and zeros of $\mathbb{E}(\lambda )$ coincide with eigenvalues
of the operator (\ref{1-3}), see the review article \cite{San} and references therein.
In \cite{Kap} and \cite{Gar}, it is proved for exponentially decaying potentials that 
$\mathbb{E}(\lambda )$ has an analytic continuation from the lower half
plane to the upper half plane through the positive real axis, whose zeros give resonance poles.

One of the difficulties when investigating properties of the Evans function is that we have to compactifying the equation
(\ref{1-4}) by attaching two points $x = \pm \infty$ because $\mu_{\pm}$ are defined by using the boundary conditions
at $x = \pm \infty$.

In the present paper, a certain function $\mathbb{D}(\lambda )$, which is holomorphic on the Riemann surface of the 
the generalized resolvent $\mathcal{R}_\lambda $, is constructed so that $\mathbb{D}(\lambda )=0$
if and only if $\lambda $ is a generalized eigenvalue.
In particular, when the Evans function is well defined, it is proved that
\begin{equation}
\mathbb{D}(\lambda ) = \frac{1}{2\sqrt{-\lambda }}\mathbb{E}(\lambda ).
\end{equation}
As a consequence, it turns out that 
the Evans function has an analytic continuation to the Riemann surface of $\mathcal{R}_\lambda $.
This gives a generalization of the results of \cite{Kap} and \cite{Gar}, in which the existence of 
the analytic continuation of $\mathbb{E}(\lambda )$ is proved only for exponentially decaying potentials.
Note also that the function $\mathbb{D}(\lambda )$ can be defined even when solutions $\mu_{\pm}$ satisfying
(\ref{1-1}), (\ref{1-2}) do not exist (this may happen when $\int_{\mathbf{R}}\! |V(x)| dx = \infty$).
Properties of $\mathbb{D}(\lambda )$ follow from those of the generalized resolvent, and we need not
consider the compactification of Eq.(\ref{1-4}).
Our results reveal that the existence of analytic continuations of the Evans functions essentially relays on the 
fact that the resolvent operator of a differential operator has an $X'$-valued 
analytic continuation if a suitable rigged Hilbert space $X \subset \mathcal{H} \subset X'$ can be constructed.

Throughout this paper, $\mathsf{D}(\cdot)$ and $\mathsf{R}(\cdot)$ denote the domain and range of an operator,
respectively.


\section{A review of the spectral theory on rigged Hilbert spaces}

\begin{figure}
\begin{center}
\includegraphics{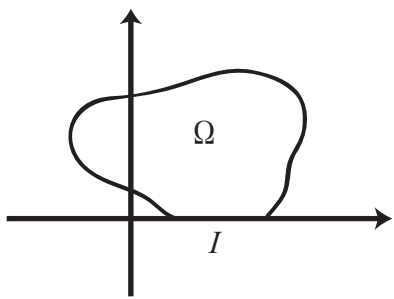}
\caption{A domain on which $E[\psi, \phi](\omega )$ is holomorphic. \label{fig1}}
\end{center}
\end{figure}

This section is devoted to a review of the spectral theory on rigged Hilbert spaces developed in \cite{Chi}.
In order to apply Schr\"{o}dinger operators, assumptions given \cite{Chi} will be slightly relaxed.
Let $\mathcal{H}$ be a Hilbert space over $\mathbf{C}$ and $H$ a selfadjoint operator densely defined on $\mathcal{H}$
with the spectral measure $\{ E(B)\}_{B\in \mathcal{B}}$; that is, $H$ is expressed as
$H = \int_{\mathbf{R}}\! \omega dE(\omega )$.
Let $K$ be some linear operator densely defined on $\mathcal{H}$.
Our purpose is to investigate spectral properties of the operator $T:=H+K$.
Let $\Omega \subset \mathbf{C}$ be a simply connected open domain in the upper half plane such that the intersection of the
real axis and the closure of $\Omega $ is a connected interval $\tilde{I}$.
Let $I = \tilde{I}\backslash \partial \tilde{I}$ be an open interval (see Fig.\ref{fig1}).
For a given $T = H+K$, we suppose that there exists a locally convex Hausdorff vector space $X(\Omega )$ 
over $\mathbf{C}$ satisfying following conditions.
\\[0.2cm]
\textbf{(X1)} $X(\Omega )$ is a dense subspace of $\mathcal{H}$.
\\
\textbf{(X2)} A topology on $X(\Omega )$ is stronger than that on $\mathcal{H}$.
\\
\textbf{(X3)} $X(\Omega )$ is a quasi-complete barreled space.
\\[0.2cm]
Let $X(\Omega )'$ be a dual space of $X(\Omega )$, the set of continuous \textit{anti}-linear functionals on $X(\Omega )$.
The paring for $(X(\Omega )', X(\Omega ))$ is denoted by $\langle \, \cdot \,|\, \cdot \, \rangle$.
For $\mu \in X(\Omega )', \phi \in X(\Omega )$ and $a\in \mathbf{C}$, we have $a \langle \mu \,|\,  \phi \rangle
=\langle a \mu \,|\,  \phi \rangle = \langle \mu \,|\,  \overline{a} \phi \rangle$.
The space $X(\Omega )'$ is equipped with the strong dual topology or the weak dual topology.
Because of (X1) and (X2), $\mathcal{H}'$, the dual of $\mathcal{H}$, is dense in $X(\Omega )'$.
Through the isomorphism $\mathcal{H} \simeq \mathcal{H}'$, we obtain the triplet
\begin{equation}
X(\Omega ) \subset \mathcal{H} \subset X(\Omega )',
\end{equation}
which is called the \textit{rigged Hilbert space} or the \textit{Gelfand triplet}.
The \textit{canonical inclusion} $i: \mathcal{H} \to X(\Omega )'$ is defined as follows; for $\psi\in \mathcal{H}$,
we denote $i(\psi)$ by $\langle \psi |$, which is defined to be
\begin{equation}
i(\psi)(\phi) = \langle \psi \,|\, \phi \rangle = (\psi, \phi),
\label{inclusion}
\end{equation}
for any $\phi \in X(\Omega )$, where $(\, \cdot\, , \, \cdot \,)$ is the inner product on $\mathcal{H}$.
The inclusion from $X(\Omega )$ into $X(\Omega )'$ is also defined as above.
Then, $i$ is injective and continuous.
The topological condition (X3) is assumed to define Pettis integrals and Taylor expansions of $X(\Omega )'$-valued
holomorphic functions.
Any complete Montel spaces, Fr\'{e}chet spaces, Banach spaces and Hilbert spaces satisfy (X3)
(we refer the reader to \cite{Tre} for basic notions of locally convex spaces, though Hilbert spaces
are mainly used in this paper).
Next, for the spectral measure $E(B)$ of $H$, we make the following analyticity conditions:
\\[0.2cm]
\textbf{(X4)} For any $\phi \in X(\Omega )$, the spectral measure $(E(B)\phi, \phi)$ is absolutely continuous on the interval $I$.
Its density function, denoted by $E[\phi, \phi](\omega )$, has an analytic continuation to $\Omega \cup I$.
\\
\textbf{(X5)} For each $\lambda \in I \cup \Omega $, the bilinear form
$E[\, \cdot\,, \, \cdot \,](\lambda ): X(\Omega ) \times X(\Omega ) \to \mathbf{C}$ is separately continuous.
\\[0.2cm]
Due to the assumption (X4) with the aid of the polarization identity, 
we can show that $(E(B)\phi, \psi)$ is absolutely continuous on $I$ for any $\phi, \psi \in X(\Omega )$.
Let $E[\phi, \psi] (\omega )$ be the density function;
\begin{equation}
d(E(\omega )\phi, \psi) = E[\phi, \psi] (\omega ) d\omega , \quad \omega \in I.
\end{equation}
Then, $E[\phi, \psi] (\omega )$ is holomorphic in $\omega \in I\cup \Omega $.
We will use the above notation for any $\omega \in \mathbf{R}$ for simplicity, although the absolute continuity is
assumed only on $I$.
Let $iX(\Omega )$ be the inclusion of $X(\Omega )$ into $X(\Omega )'$.
Define the operator $A(\lambda ) : iX(\Omega ) \to X(\Omega )'$ to be
\begin{equation}
\langle A(\lambda )\psi \,|\, \phi \rangle = \left\{ \begin{array}{ll}
\displaystyle 
     \int_{\mathbf{R}}\! \frac{1}{\lambda -\omega } E [\psi, \phi](\omega ) d\omega 
           + 2\pi \sqrt{-1} E [\psi, \phi](\lambda )  & (\lambda \in \Omega ), \\[0.4cm]
\displaystyle  \lim_{y\to -0} 
  \int_{\mathbf{R}}\! \frac{1}{x + \sqrt{-1}y -\omega } E [\psi, \phi](\omega ) d\omega  & (\lambda =x\in I), \\[0.4cm]
\displaystyle  \int_{\mathbf{R}}\! \frac{1}{\lambda -\omega } E [\psi, \phi](\omega ) d\omega
 & (\mathrm{Im}(\lambda ) < 0),
\end{array} \right.
\label{A}
\end{equation}
for any $\psi \in iX(\Omega )$ and $\phi \in X(\Omega )$.
It is known that $\langle A(\lambda )\psi \,|\, \phi \rangle$ is holomorphic on the region
$\{ \mathrm{Im}(\lambda ) < 0\} \cup \Omega \cup I$.
It is proved in \cite{Chi} that $A(\lambda ) \circ i : X(\Omega ) \to X(\Omega )'$ is continuous when
$X(\Omega )'$ is equipped with the weak dual topology.
When $\mathrm{Im}(\lambda ) < 0$, we have $\langle A(\lambda )\psi \,|\, \phi \rangle = ( (\lambda - H)^{-1}\psi, \phi )$.
In this sense, the operator $A(\lambda )$ is called the analytic continuation of the resolvent
$(\lambda -H )^{-1}$ in the generalized sense.
The operator $A(\lambda )$ plays a central role for our theory. 

Let $Q$ be a linear operator densely defined on $X(\Omega )$.
Then, the dual operator $Q'$ is defined as follows:
the domain $\mathsf{D}(Q')$ of $Q'$ is the set of elements $\mu\in X(\Omega )'$ such that the mapping
$\phi \mapsto \langle \mu \,|\, Q\phi \rangle$ from $X(\Omega )$ into $\mathbf{C}$ is continuous.
Then, $Q' : \mathsf{D}(Q') \to X(\Omega )'$ is defined by 
$\langle Q' \mu \,|\, \phi \rangle = \langle \mu \,|\, Q\phi \rangle$.
The (Hilbert) adjoint $Q^*$ of $Q$ is defined through $(Q\phi, \psi) = (\phi , Q^* \psi)$ as usual
when $Q$ is densely defined on $\mathcal{H}$.
If $Q^*$ is densely defined on $X(\Omega )$, its dual $(Q^*)'$ is well defined, which is denoted by $Q^\times$.
Then, $Q^\times = (Q^*)'$ is an extension of $Q$ which satisfies $i\circ Q = Q^\times \circ i\,|_{\mathsf{D}(Q)}$.
For the operators $H$ and $K$, we suppose that
\\[0.2cm]
\textbf{(X6)} there exists a dense subspace $Y$ of $X(\Omega )$ such that $HY \subset X(\Omega )$.
\\
\textbf{(X7)} $K$ is $H$-bounded and there exists a dense subspace $Y$ of $X(\Omega )$ such that $K^*Y \subset X(\Omega )$.
\\
\textbf{(X8)} $K^\times A(\lambda ) iX(\Omega ) \subset iX(\Omega )$ for any $\lambda \in \{ \mathrm{Im}(\lambda )<0\}
\cup I \cup \Omega $.
\\[0.2cm]
Due to (X6) and (X7), we can show that $H^\times$, $K^\times$ and $T^\times$ are densely defined on $X(\Omega )'$.
In particular, $\mathsf{D}(H^\times) \supset iY, \mathsf{D}(K^\times) \supset iY$ and
$\mathsf{D}(T^\times) \supset iY$.
When $H$ and $K$ are continuous on $X(\Omega )$, (X6) and (X7) are satisfied with $Y = X(\Omega )$.
Then, $H^\times $ and $T^\times$ are continuous on $X(\Omega )'$.
Recall that $K$ is called $H$-bounded if $K(\lambda -H)^{-1}$ is bounded on $\mathcal{H}$.
In particular, $K(\lambda -H)^{-1}\mathcal{H}\subset \mathcal{H}$.
Since $A(\lambda )$ is the analytic continuation of $(\lambda -H)^{-1}$ as an operator from $iX(\Omega )$,
(X8) gives an ``analytic continuation version" of the assumption that $K$ is $H$-bounded.
In \cite{Chi}, the spectral theory of the operator $T=H+K$ is developed under the assumptions (X1) to (X8).
However, we will show that a Schr\"{o}dinger operator with a dilation analytic potential does not satisfy the assumption (X8).
Thus, we make the following condition instead of (X8).
In what follows, put $\hat{\Omega } = \Omega \cup I \cup \{ \lambda \, | \, \mathrm{Im}(\lambda )<0\}$.

Suppose that there exists a locally convex Hausdorff vector space $Z(\Omega )$ satisfying following conditions:
\\
\textbf{(Z1)} $X(\Omega )$ is a dense subspace of $Z(\Omega )$ and 
the topology of $X(\Omega )$ is stronger than that of $Z(\Omega )$.
\\
\textbf{(Z2)} $Z(\Omega )$ is a quasi-complete barreled space.
\\
\textbf{(Z3)} The canonical inclusion $i : X(\Omega ) \to X(\Omega )'$ is continuously extended
to a mapping $j : Z(\Omega ) \to X(\Omega )'$.
\\
\textbf{(Z4)} For any $\lambda \in \hat{\Omega }$, the operator $A(\lambda ) : iX(\Omega ) \to X(\Omega )'$
is extended to an operator from $jZ(\Omega )$ into $X(\Omega )'$ so that $A(\lambda ) \circ j : Z(\Omega ) \to X(\Omega )'$
is continuous if $X(\Omega )'$ is equipped with the weak dual topology.
\\
\textbf{(Z5)} For any $\lambda \in \hat{\Omega }$, $K^\times A(\lambda )j Z(\Omega ) \subset jZ(\Omega )$
and $j^{-1}K^\times A(\lambda )j$ is continuous on $Z(\Omega )$.
\begin{eqnarray*}
\begin{array}{ccccc}
X(\Omega ) & \subset & \mathcal{H} & \subset & X(\Omega )'\\
\rotatebox{90}{$\supset$} & & & & \rotatebox{90}{$\subset$} \\
Z(\Omega ) & &\longrightarrow   & & jZ(\Omega )\\
\end{array}
\end{eqnarray*}

If $Z(\Omega ) = X(\Omega )$, then (Z1) to (Z5) are reduced to (X1) to (X8), 
and the results obtained in \cite{Chi} are recovered.
If $X(\Omega ) \subset Z(\Omega ) \subset \mathcal{H}$,
then $X(\Omega )$ does not play a role; we should use the triplet
$Z(\Omega ) \subset \mathcal{H} \subset Z(\Omega )'$ from the beginning.
Thus we are interested in the situation $Z(\Omega ) \nsubset \mathcal{H}$.
In what follows, the extension $j$ of $i$ is also denoted by $i$ for simplicity.
Let us show the same results as \cite{Chi} under the assumptions (X1) to (X7) and (Z1) to (Z5).
\\[0.2cm]
\textbf{Lemma \thedef.}
\\
(i) For each $\phi \in Z(\Omega )$, $A(\lambda )i\phi$ is an $X(\Omega )'$-valued holomorphic function
in $\lambda \in \hat{\Omega }$.
\\
(ii) Define the operators $A^{(n)}(\lambda ) : i X(\Omega ) \to X(\Omega )'$ to be
\begin{eqnarray}
\langle A^{(n)}(\lambda )\psi \,|\, \phi \rangle = \left\{ \begin{array}{l}
\displaystyle 
     \int_{\mathbf{R}}\! \frac{1}{(\lambda -\omega )^n} E [\psi, \phi](\omega ) d\omega 
    + 2\pi \sqrt{-1} \frac{(-1)^{n-1}}{(n-1)!} \frac{d^{n-1}}{dz^{n-1}}\Bigl|_{z=\lambda }
               E [\psi, \phi](z), \,\, (\lambda \in \Omega ), \\[0.4cm]
\displaystyle  \lim_{y\to -0} 
  \int_{\mathbf{R}}\! \frac{1}{(x + \sqrt{-1}y -\omega )^n} E [\psi, \phi](\omega ) d\omega, \quad (\lambda =x\in I), \\[0.4cm]
\displaystyle  \int_{\mathbf{R}}\! \frac{1}{(\lambda -\omega )^n} E [\psi, \phi](\omega ) d\omega,
\quad (\mathrm{Im}(\lambda ) < 0),
\end{array} \right.
\label{An}
\end{eqnarray}
for $n=1, 2\cdots $.
Then, $A^{(n)}(\lambda ) \circ i$ has a continuous extension $A^{(n)}(\lambda ) \circ i : Z(\Omega ) \to X(\Omega )'$,
and $A(\lambda )i\phi$ is expanded in a Taylor series as
\begin{equation}
A(\lambda )i\phi = \sum^\infty_{j=0}(\lambda _0 - \lambda )^j A^{(j+1)}(\lambda _0)i\phi, \quad \phi \in Z(\Omega ),
\label{Taylor}
\end{equation}
which converges with respect to the strong dual topology on $X(\Omega )'$.
\\
(iii) When $\mathrm{Im} (\lambda ) < 0$, $A(\lambda )\circ i\phi = i\circ (\lambda -H)^{-1} \phi$ 
for $\phi \in X(\Omega )$.
\\[0.2cm]
\textbf{Proof.}
(i) In \cite{Chi}, $\langle A(\lambda )i\phi \,|\, \psi \rangle$ is proved to be holomorphic in $\lambda \in \hat{\Omega }$
for any $\phi, \psi \in X(\Omega )$.
Since $X(\Omega )$ is dense in $Z(\Omega )$, Montel theorem proves that $\langle A(\lambda )i\phi \,|\, \psi \rangle$
is holomorphic for $\phi\in Z(\Omega )$ and $\psi \in X(\Omega )$.
This implies that $A(\lambda )i\phi$ is a weakly holomorphic $X(\Omega )'$-valued function.
Since $X(\Omega )$ is barreled, Thm.A.3 of \cite{Chi} concludes that $A(\lambda )i\phi$ is strongly holomorphic.
(ii) In \cite{Chi}, Eq.(\ref{Taylor}) is proved for $\phi \in X(\Omega )$.
Again Montel theorem is applied to show the same equality for $\phi \in Z(\Omega )$.
(iii) This follows from the definition of $A(\lambda )$. \hfill $\blacksquare$
\\[-0.2cm]

Lemma 2.1 means that $A(\lambda )$ gives an analytic continuation of the resolvent $(\lambda -H)^{-1}$ from
the lower half plane to $\Omega $ as an $X(\Omega )'$-valued function.
Similarly, $A^{(n)}(\lambda )$ is an analytic continuation of $(\lambda -H)^{-n}$.
$A^{(1)}(\lambda )$ is also denoted by $A(\lambda )$ as before.
Next, let us define an analytic continuation of the resolvent of $T=H+K$.
Due to (Z5), $id - K^\times A(\lambda )$ is an operator on $iZ(\Omega )$.
It is easy to verify that $id - K^\times A(\lambda )$ is injective if and only if 
$id - A(\lambda )K^\times$ is injective on $\mathsf{R}(A(\lambda )) = A(\lambda )iZ(\Omega )$.
\\[0.2cm]
\textbf{Definition \thedef.}
If the inverse $(id - K^\times A(\lambda ))^{-1}$ exists on $iZ(\Omega )$,
define the generalized resolvent $\mathcal{R}_\lambda : iZ(\Omega ) \to X(\Omega )'$ of $T$ to be
\begin{equation}
\mathcal{R}_\lambda =A(\lambda ) \circ (id - K^\times A(\lambda ))^{-1}
      = (id - A(\lambda )K^\times)^{-1} \circ A(\lambda ), \quad \lambda \in \hat{\Omega }.
\label{resolvent}
\end{equation}
Although $\mathcal{R}_\lambda $ is not a continuous operator in general,
the composition $\mathcal{R}_\lambda \circ i : Z(\Omega ) \to X(\Omega )'$ may be continuous:
\\[0.2cm]
\textbf{Definition \thedef.}
The generalized resolvent set $\hat{\varrho}(T)$ is defined to be the set of points $\lambda \in \hat{\Omega }$ satisfying following:
there is a neighborhood $V_\lambda \subset \hat{\Omega }$ of $\lambda $ such that for any $\lambda ' \in V_\lambda $,
$\mathcal{R}_{\lambda '} \circ i$ is a densely 
defined continuous operator from $Z(\Omega )$ into $X(\Omega )'$, where $X(\Omega )'$ is equipped with the weak dual topology,
and the set $\{ \mathcal{R}_{\lambda '} \circ i(\psi)\}_{\lambda ' \in V_\lambda }$ is bounded in $X(\Omega )'$ for each $\psi \in Z(\Omega )$.
The set $\hat{\sigma }(T):= \hat{\Omega }\backslash \hat{\varrho}(T)$ is called the \textit{generalized spectrum} of $T$.
The \textit{generalized point spectrum} $\hat{\sigma }_p(T)$ is the set of points $\lambda \in \hat{\sigma }(T)$ at which
$id - K^\times A(\lambda )$ is not injective.
The \textit{generalized residual spectrum} $\hat{\sigma }_r(T)$ is the set of points $\lambda \in \hat{\sigma }(T)$ 
such that the domain of $\mathcal{R}_\lambda \circ i$ is not dense in $Z(\Omega )$.
The \textit{generalized continuous spectrum} is defined to be 
$\hat{\sigma }_c(T) = \hat{\sigma }(T)\backslash (\hat{\sigma }_p(T)\cup \hat{\sigma }_r(T))$.
\\

We can show that if $Z(\Omega )$ is a Banach space, $\lambda \in \hat{\varrho}(T)$ if and only if 
$id - i^{-1}K^\times A(\lambda )i$ has a continuous inverse on $Z(\Omega )$ (Prop.3.18 of \cite{Chi}).
The next theorem is proved in the same way as Thm.3.12 of \cite{Chi}.
\\[0.2cm]
\textbf{Theorem \thedef \,\cite{Chi}.} Suppose (X1) to (X7) and (Z1) to (Z5).
\\
(i) For each $\phi \in Z(\Omega )$, $\mathcal{R}_\lambda i\phi $ is an $X(\Omega )'$-valued holomorphic function
in $\lambda \in \hat{\varrho} (T)$.
\\
(ii) Suppose $\mathrm{Im} (\lambda ) < 0$, $\lambda \in \hat{\varrho} (T)$ and $\lambda \in \varrho (T)$
(the resolvent set of $T$ in $\mathcal{H}$-sense).
Then, $\mathcal{R}_\lambda \circ i\phi= i\circ (\lambda -T)^{-1} \phi$ for any $\phi \in X(\Omega )$.
In particular, $\langle \mathcal{R}_\lambda i\phi \,|\, \psi \rangle$ is an analytic continuation of 
$((\lambda -T)^{-1}\phi, \psi)$ from the lower half plane to $I\cup \Omega $ for any $\phi, \psi \in X(\Omega )$.
\\

Next, we define the operator $B^{(n)}(\lambda ) : \mathsf{D}(B^{(n)}(\lambda ))
\subset X(\Omega )' \to X(\Omega )'$ to be
\begin{equation}
B^{(n)}(\lambda )=id - A^{(n)}(\lambda )K^\times (\lambda -H^\times)^{n-1}.
\label{2-7}
\end{equation}
The domain $\mathsf{D}(B^{(n)}(\lambda ))$ is the set of $\mu \in X(\Omega )'$ such that
$K^\times (\lambda -H^\times)^{n-1} \mu \in iZ(\Omega )$.
\\[0.2cm]
\textbf{Definition \thedef.} 
A point $\lambda$ in $\hat{\sigma }_p(T)$ is called a generalized eigenvalue (resonance pole) of the operator $T$.
The generalized eigenspace of $\lambda $ is defined by
\begin{equation}
V_\lambda = \bigcup_{m\geq 1}\mathrm{Ker}\, B^{(m)}(\lambda ) \circ B^{(m-1)}(\lambda ) \circ \cdots \circ B^{(1)}(\lambda ).
\label{space}
\end{equation}
We call $\mathrm{dim} V_\lambda $ the algebraic multiplicity of the generalized eigenvalue $\lambda $.
In particular, a nonzero solution $\mu \in X(\Omega )'$ of the equation
\begin{equation}
B^{(1)}(\lambda )\mu = (id - A(\lambda )K^\times) \mu = 0
\end{equation}
is called a \textit{generalized eigenfunction} associated with the generalized eigenvalue $\lambda $.
\\[0.2cm]
\textbf{Theorem \thedef \,\cite{Chi}.} For any $\mu \in V_\lambda $, there exists an integer $M$
such that $(\lambda -T^\times) ^M \mu = 0$.
In particular, a generalized eigenfunction $\mu$ satisfies $(\lambda -T^\times) \mu = 0$.
\\

This implies that $\lambda $ is indeed an eigenvalue of the dual operator $T^\times$.
In general, $\hat{\sigma }_p(T)$ is a proper subset of $\sigma _p(T^\times)$ (the set of eigenvalues of $T^\times$),
and $V_\lambda $ is a proper subspace of the eigenspace $\bigcup_{m\geq 1} \mathrm{Ker} (\lambda -T^\times)^m$ of $T^\times$.

Let $\Sigma \subset \hat{\sigma}(T)$ be a bounded subset of the generalized spectrum,
which is separated from the rest of the spectrum by a simple closed curve $\gamma \subset \hat{\Omega}$.
Define the operator $\Pi_\Sigma : iZ(\Omega ) \to X(\Omega )'$ to be
\begin{equation}
\Pi_\Sigma \phi= \frac{1}{2\pi \sqrt{-1}} \int_{\gamma }\! \mathcal{R}_\lambda \phi \,d\lambda,
\quad \phi \in iZ(\Omega ),
\label{2-10}
\end{equation}
which is called the \textit{generalized Riesz projection} for $\Sigma $.
The integral in the right hand side is well defined as the Pettis integral.
We can show that $\Pi_\Sigma \circ i$ is a continuous operator from $Z(\Omega )$ into $X(\Omega )'$
equipped with the weak dual topology.
Note that $\Pi_\Sigma \circ \Pi_\Sigma  = \Pi_\Sigma $ does not hold because the composition $\Pi_\Sigma \circ \Pi_\Sigma $
is not defined.
Nevertheless, we call it the projection because it is proved in 
Prop.3.14 of \cite{Chi} that $\Pi_\Sigma (iZ(\Omega )) \cap (id - \Pi_\Sigma)(iZ(\Omega )) = \{ 0\}$
and the direct sum satisfies
\begin{equation}
iZ(\Omega ) \subset \Pi_\Sigma (iZ(\Omega )) \oplus (id - \Pi_\Sigma)(iZ(\Omega )) \subset X(\Omega )'.
\end{equation}
Let $\lambda _0$ be an isolated generalized eigenvalue, which is separated from the rest of the generalized 
spectrum by a simple closed curve $\gamma _0 \subset \hat{\Omega}$.
Let
\begin{equation}
\Pi_0 = \frac{1}{2\pi \sqrt{-1}} \int_{\gamma _0}\! \mathcal{R}_\lambda d\lambda ,
\label{pro}
\end{equation}
be a projection for $\lambda _0$ and $V_0 = \bigcup_{m\geq 1} \mathrm{Ker}\, B^{(m)}(\lambda _0) \circ
\cdots \circ B^{(1)}(\lambda _0)$ a generalized eigenspace of $\lambda _0$.
\\[0.2cm]
\textbf{Theorem \thedef \,\cite{Chi}.}
If $\Pi_0 iZ(\Omega ) = \mathsf{R}(\Pi_0)$ is finite dimensional, then $\Pi_0 iZ(\Omega ) = V_0$.
\\[0.2cm]
Note that $\Pi_0 iZ(\Omega ) = \Pi_0 iX(\Omega )$ when $\Pi_0 iZ(\Omega )$ is finite dimensional because
$X(\Omega )$ is dense in $Z(\Omega )$.
Then, the above theorem is proved in the same way as the proof of Thm.3.16 of \cite{Chi}.
\\[0.2cm]
\textbf{Theorem \thedef \,\cite{Chi}.}
In addition to (X1) to (X7) and (Z1) to (Z5), suppose that
\\[0.2cm]
\textbf{(Z6)} $i^{-1} K^\times A(\lambda )i : Z(\Omega ) \to Z(\Omega )$ 
is a compact operator uniformly in $\lambda \in \hat{\Omega }$.
\\[0.2cm]
Then, the following statements are true.
\\ 
(i) For any compact set $D \subset \hat{\Omega }$, the number of generalized eigenvalues in $D$ is finite
(thus $\hat{\sigma }_p(T)$ consists of a countable number of generalized eigenvalues and they may accumulate
only on the boundary of $\hat{\Omega }$ or infinity).
\\
(ii) For each $\lambda _0\in \hat{\sigma }_p(T)$, the generalized eigenspace $V_0$ is of finite dimensional and 
$\Pi_0iZ(\Omega ) = V_0$.
\\
(iii) $\hat{\sigma }_c(T) = \hat{\sigma }_r(T) = \emptyset$.

Recall that a linear operator $L$ from a topological vector space $X_1$ to another topological vector space $X_2$
is said to be compact if there exists a neighborhood $U\subset X_1$ such that $LU \subset X_2$ is relatively compact.
When $L=L(\lambda )$ is parameterized by $\lambda $, it is said to be compact uniformly in $\lambda $
if such a neighborhood $U$ is independent of $\lambda $.
When the domain $X_1$ is a Banach space, $L(\lambda )$ is compact uniformly in $\lambda $
if and only if $L(\lambda )$ is compact for each $\lambda $.
The above theorem is also proved in a similar manner to the proof of Thm.3.19 of \cite{Chi}.
It is remarkable that $\hat{\sigma }_c(T) =\emptyset$ even if $T$ has the continuous spectrum in $\mathcal{H}$-sense.

When we emphasize the choice of $Z(\Omega )$, $\hat{\sigma }(T)$ is also denoted by $\hat{\sigma }(T; Z(\Omega ))$.
Now suppose that two vector spaces $Z_1(\Omega )$ and $Z_2(\Omega )$ satisfy the 
assumptions (Z1) to (Z5) with a common $X(\Omega )$.
Then, two generalized spectra $\hat{\sigma }(T;Z_1(\Omega ))$ and $\hat{\sigma }(T;Z_2(\Omega ))$ for 
$Z_1(\Omega )$ and $Z_2(\Omega )$ are defined, respectively.
Let us consider the relationship between them.
\\[0.2cm]
\textbf{Proposition \thedef.}
Suppose that $Z_2(\Omega )$ is a dense subspace of $Z_1(\Omega )$ and the topology on $Z_2(\Omega )$ is stronger than
that on $Z_1(\Omega )$. Then, the following holds.
\\
(i) $\hat{\sigma }(T;Z_2(\Omega )) \subset \hat{\sigma }(T;Z_1(\Omega ))$.
\\
(ii) Let $\Sigma $ be a bounded subset of $\hat{\sigma} (T; Z_1(\Omega ))$ which is separated from the rest of 
the spectrum by a simple closed curve $\gamma $.
Then, there exists a point of $\hat{\sigma }(T; Z_2(\Omega ))$ inside $\gamma $.
In particular, if $\lambda $ is an isolated point of $\hat{\sigma} (T; Z_1(\Omega ))$,
then $\lambda \in \hat{\sigma }(T; Z_2(\Omega ))$.
\\[0.2cm]
\textbf{Proof.}
(i) Suppose that $\lambda \notin \hat{\sigma }(T; Z_1(\Omega ))$.
Then, there is a neighborhood $V_\lambda$ of $\lambda $ such that $\mathcal{R}_{\lambda '} \circ i$ is a continuous operator 
from $Z_1(\Omega )$ into $X(\Omega )'$ for any $\lambda ' \in V_\lambda $, and the set 
$\{ \mathcal{R}_{\lambda '} \circ i \psi \}_{\lambda ' \in V_\lambda }$ is bounded in $X(\Omega )'$ for each $\psi \in Z_1(\Omega )$.
Since the topology on $Z_2(\Omega )$ is stronger than that on $Z_1(\Omega )$, $\mathcal{R}_{\lambda '} \circ i$ is a continuous operator 
from $Z_2(\Omega )$ into $X(\Omega )'$ for any $\lambda ' \in V_\lambda $, and the set 
$\{ \mathcal{R}_{\lambda '} \circ i \psi \}_{\lambda ' \in V_\lambda }$ is bounded in $X(\Omega )'$ for each $\psi \in Z_2(\Omega )$.
This proves that $\lambda \notin \hat{\sigma }(T; Z_2(\Omega ))$.

(ii) Let $\Pi_\Sigma$ be the generalized Riesz projection for $\Sigma$.
Since $\Sigma \subset \hat{\sigma }(T; Z_1(\Omega ))$, $\Pi_{\Sigma} iZ_1(\Omega )\neq \{ 0 \}$.
Then, $\Pi_{\Sigma} iZ_2(\Omega )\neq \{ 0\}$ because $Z_2(\Omega )$ is dense in $Z_1(\Omega )$.
This shows that the closed curve $\gamma $ encloses a point of $\hat{\sigma }(T; Z_2(\Omega ))$. \hfill $\blacksquare$


\section{An application to Schr\"{o}dinger operators}

In this section, we consider a Schr\"{o}dinger operator of the form $T = -\Delta + V$,
where $\Delta$ is the Laplace operator on $\mathbf{R}^m$ defined to be
\begin{equation}
\Delta = \frac{\partial ^2}{\partial x_1^2} + \frac{\partial ^2}{\partial x_2^2}
 + \cdots  + \frac{\partial ^2}{\partial x_m^2},
\end{equation}
and $V$ is the multiplication operator by a function $V : \mathbf{R}^m \to \mathbf{C}$;
\begin{equation}
(V\phi) (x) = V(x) \cdot \phi (x), \quad x\in \mathbf{R}^m.
\end{equation}
Put $\mathcal{H} = L^2 (\mathbf{R}^m)$.
The domain of $\Delta$ is the Sobolev space $H^2 (\mathbf{R}^m)$.
Then, $\Delta$ is a selfadjoint operator densely defined on $\mathcal{H}$.
In what follows, we denote $-\Delta$ and $V$ by $H$ and $K$, respectively.
Our purpose is to investigate an operator $T = H+K$ with suitable assumptions for $K=V$.

Define the Fourier transform and the inverse Fourier transform to be
\begin{equation}
\mathcal{F}[u](\xi) = \frac{1}{(2\pi)^{m/2}} \int_{\mathbf{R}^m}\! u(x)e^{-\sqrt{-1}x\cdot \xi} dx, \quad 
\mathcal{F}^{-1}[\hat{u}](x) = \frac{1}{(2\pi)^{m/2}} \int_{\mathbf{R}^m}\! \hat{u}(\xi)e^{\sqrt{-1}x\cdot \xi} d\xi,
\end{equation}
where $x = (x_1, \cdots  ,x_m)$ and $\xi = (\xi_1, \cdots , \xi_m)$.
The resolvent of $H$ is given by
\begin{eqnarray*}
(\lambda -H)^{-1}\psi (x) &=&
 \frac{1}{(2\pi)^{m/2}} \int_{\mathbf{R}^m}\! \frac{1}{\lambda -|\xi|^2} e^{\sqrt{-1}x\cdot \xi} \mathcal{F}[\psi ](\xi) d\xi,
\end{eqnarray*}
where $|\xi|^2 = \xi_1^2 + \cdots + \xi^2_m$.
Let $S^{m-1} \subset \mathbf{R}^m$ be the $(m-1)$-dimensional unit sphere.
For a point $\xi \in \mathbf{R}^m$, put $\xi = r\omega $ with $r\geq 0,\, \omega \in S^{m-1}$.
Then, $(\lambda -H)^{-1}\psi  (x)$ is rewritten as
\begin{eqnarray}
(\lambda -H)^{-1}\psi  (x) &=&
 \frac{1}{(2\pi)^{m/2}}\int_{S^{m-1}}\! d\omega 
\int^\infty_{0}\! \frac{1}{\lambda -r^2} e^{\sqrt{-1}r x\cdot \omega } \mathcal{F}[\psi ](r\omega ) r^{m-1}dr \\
&=&  \frac{1}{(2\pi)^{m/2}} 
\int^{\infty}_{0}\! \frac{1}{\lambda -r}\left( \int_{S^{m-1}}\! \frac{\sqrt{r}^{m-2}}{2} e^{\sqrt{-1}\sqrt{r} x\cdot \omega}
\mathcal{F}[\psi ] (\sqrt{r}\omega ) d\omega  \right)  dr,
\label{5-7}
\end{eqnarray}
which gives the spectral representation of the resolvent.
This is an $L^2(\mathbf{R}^m)$-valued holomorphic function in $\{ \lambda \, | \, -2\pi < \mathrm{arg} (\lambda ) < 0\}$
for each $\psi \in L^2(\mathbf{R}^m)$.
The positive real axis $\mathrm{arg} (\lambda ) = 0$ is the essential spectrum of $H$.
Let $\Omega $ be an open domain on the upper half plane as in Fig.\ref{fig1}.
If the function $f(z):= \mathcal{F}[\psi ] (\sqrt{z}\omega )$ is holomorphic on $\Omega $,
then the above quantity has an analytic continuation with respect to $\lambda $ 
from the sector $-\varepsilon < \mathrm{arg} (\lambda ) < 0$ on the lower half plane to $\Omega $ as
\begin{eqnarray*}
& &  \frac{1}{(2\pi)^{m/2}}\int^{\infty}_{0}\! \frac{1}{\lambda -r}\left( \int_{S^{m-1}}\! \frac{\sqrt{r}^{m-2}}{2} e^{\sqrt{-1}\sqrt{r} x\cdot \omega}
\mathcal{F}[\psi ] (\sqrt{r}\omega ) d\omega  \right)  dr \\
& & \quad \quad   + \frac{\pi \sqrt{-1}}{(2\pi)^{m/2}}\sqrt{\lambda }^{m-2}\!\!\int_{S^{m-1}}\! \!
 e^{\sqrt{-1}\sqrt{\lambda }x \cdot \omega } \mathcal{F}[\psi ] (\sqrt{\lambda }\omega ) d\omega ,
\end{eqnarray*}
which is not included in $L^2(\mathbf{R}^m)$ in general.
Suppose for simplicity that $\mathcal{F}[\psi ] (\sqrt{z }\omega )$ is an entire function with respect 
to $\sqrt{z}$; that is, $f(z)= \mathcal{F}[\psi ] (\sqrt{z}\omega )$ is holomorphic on the Riemann surface of $\sqrt{z}$.
Then, the above function exists for $\{ \lambda \, | \, 0 \leq \mathrm{arg} (\lambda ) < 2\pi \}$.
Furthermore, the analytic continuation of the above function from the sector 
$2 \pi -\varepsilon < \mathrm{arg} (\lambda ) < 2 \pi$ to the 
upper half plane through the ray $\mathrm{arg} (\lambda ) = 2\pi$ is given by
\begin{eqnarray*}
 \frac{1}{(2\pi)^{m/2}}\int^{\infty}_{0}\! \frac{1}{\lambda -r}\left( \int_{S^{m-1}}\! \frac{\sqrt{r}^{m-2}}{2} e^{\sqrt{-1}\sqrt{r} x\cdot \omega}
\mathcal{F}[\psi ] (\sqrt{r}\omega ) d\omega  \right)  dr \quad ( = (\lambda -H)^{-1}\psi  (x) ),
\end{eqnarray*}
when $m$ is an odd integer, and given by
\begin{eqnarray*}
& &  \frac{1}{(2\pi)^{m/2}}\int^{\infty}_{0}\! \frac{1}{\lambda -r}\left( \int_{S^{m-1}}\! \frac{\sqrt{r}^{m-2}}{2} e^{\sqrt{-1}\sqrt{r} x\cdot \omega}
\mathcal{F}[\psi ] (\sqrt{r}\omega ) d\omega  \right)  dr \\
& & \quad \quad   + \frac{2\pi \sqrt{-1}}{(2\pi)^{m/2}}\sqrt{\lambda }^{m-2}\!\!\int_{S^{m-1}}\! \!
 e^{\sqrt{-1}\sqrt{\lambda }x \cdot \omega } \mathcal{F}[\psi ] (\sqrt{\lambda }\omega ) d\omega ,
\end{eqnarray*}
when $m$ is an even integer.
Repeating this procedure shows that the analytic continuation $A(\lambda )i (\psi)$ of $(\lambda -H)^{-1}\psi$
in the generalized sense is given by
\begin{equation}
A(\lambda )i (\psi) (x) = 
\left\{ \begin{array}{ll}
(\lambda -H)^{-1}\psi (x) & (-2\pi < \mathrm{arg} (\lambda ) < 0), \\[0.2cm]
\displaystyle (\lambda -H)^{-1}\psi (x) & \\
\displaystyle + 
\frac{\pi \sqrt{-1}}{(2\pi)^{m/2}}\sqrt{\lambda }^{m-2}\!\!\int_{S^{m-1}}\! \!
 e^{\sqrt{-1}\sqrt{\lambda }x \cdot \omega } \mathcal{F}[\psi] (\sqrt{\lambda }\omega ) d\omega
 & (0 < \mathrm{arg} (\lambda ) < 2\pi ),\\
\end{array} \right.
\label{3-9}
\end{equation}
which is defined on the Riemann surface of $\sqrt{\lambda }$, when $m$ is an odd integer,
and given by 
\begin{equation}
A(\lambda )i (\psi) (x) = 
(\lambda -H)^{-1}\psi (x)+ 
\frac{n \cdot \pi \sqrt{-1}}{(2\pi)^{m/2}}\sqrt{\lambda }^{m-2}\!\!\int_{S^{m-1}}\! \!
 e^{\sqrt{-1}\sqrt{\lambda }x \cdot \omega } \mathcal{F}[\psi] (\sqrt{\lambda }\omega ) d\omega,
\end{equation}
for $2\pi (n-1) < \mathrm{arg} (\lambda ) < 2\pi n$, which is defined on the logarithmic Riemann surface,
when $m$ is an even integer.
In what follows, the plane $P_1=\{ \lambda \, | \, -2\pi < \mathrm{arg} (\lambda ) < 0 \}$ is 
referred to as the first Riemann sheet, on which $A(\lambda )i (\psi) (x)$ coincides with the 
resolvent $(\lambda -H)^{-1}\psi (x)$ in $L^2(\mathbf{R}^m)$-sense.
The plane $P_n=\{ \lambda \, | \, 2\pi (n-2) < \mathrm{arg} (\lambda ) < 2\pi (n-1) \}$ is 
referred to as the $n$-th Riemann sheet, on which $A(\lambda )i (\psi) (x)$ is not included in $L^2(\mathbf{R}^m)$.

Once the operator $K$ is given, we should find spaces $X(\Omega )$ and $Z(\Omega )$ so that the assumptions (X1) to (X7)
and (Z1) to (Z5) are satisfied.
Then, the above $A(\lambda )i(\psi)$ can be regarded as an $X(\Omega )'$-valued holomorphic function.
In this paper, we will give two examples.
In Sec.3.1, we consider a potential $V(x)$ which decays exponentially as $|x| \to \infty$.
In this case, we can find a space $X(\Omega )$ satisfying (X1) to (X8).
Thus we need not introduce a space $Z(\Omega )$.
In Sec.3.2, a dilation analytic potential is considered.
In this case, (X8) is not satisfied and we have to find a space $Z(\Omega )$ satisfying (Z1) to (Z5).
For an exponentially decaying potential, the formulation using a rigged Hilbert space is well known for experts.
For a dilation analytic potential, our formulation based on a rigged Hilbert space is new;
in the literature, a resonance pole for such a potential is treated by using the spectral deformation technique \cite{His}.
In our method, we need not introduce any spectral deformations.


\subsection{Exponentially decaying potentials}

Let $a>0$ be a positive number.
For the function $V$, we suppose that
\begin{equation}
e^{2a|x|} V(x) \in L^2(\mathbf{R}^m).
\label{4-11}
\end{equation}
Thus $V(x)$ has to decay with the exponential rate. 
For this $a>0$, let $X(\Omega ):=L^2(\mathbf{R}^m, e^{2a|x|}dx)$ be the weighted Lebesgue space.
It is known that the dual space $X(\Omega )'$ equipped with the strong dual topology is 
identified with the weighted Lebesgue space $L^2(\mathbf{R}^m, e^{-2a|x|}dx)$.
Let us show that the rigged Hilbert space 
\begin{equation}
L^2(\mathbf{R}^m, e^{2a|x|}dx) \subset L^2(\mathbf{R}^m) \subset L^2(\mathbf{R}^m, e^{-2a|x|}dx)
\end{equation}
satisfies the assumptions (X1) to (X8).
In what follows, we suppose that $m$ is an odd integer for simplicity.
The even integer case is treated in the same way.

\begin{figure}
\begin{center}
\includegraphics{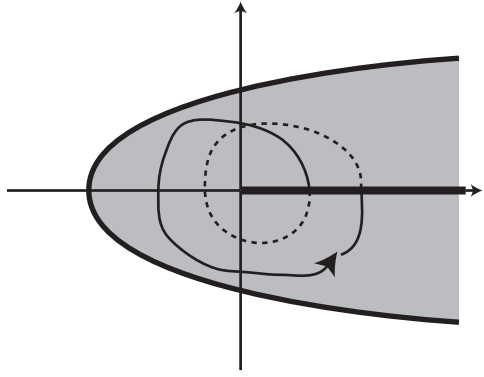}
\caption{The Riemann surface of the generalized resolvent of the Schr\"{o}dinger operator
on an odd dimensional space with an exponentially decaying potential.
The origin is a branch point of $\sqrt{z}$.
The continuous spectrum in $L^2(\mathbf{R}^m)$-sense is regarded as the branch cut.
On the region $\{ \lambda \, | \, -2\pi < \mathrm{arg} (\lambda ) < 0\}$,
the generalized resolvent $A(\lambda )$ coincides with the usual resolvent in $L^2(\mathbf{R}^m)$-sense,
while on the region $\{ \lambda \, | \, 0 < \mathrm{arg} (\lambda ) < 2\pi\}$,
the generalized resolvent is given by the second line of Eq.(\ref{3-12}).
 \label{fig2}}
\end{center}
\end{figure}

It is known that for any $\psi\in L^2(\mathbf{R}^m, e^{2a|x|}dx)$, the function $\mathcal{F}[\psi](r\omega )$ has an 
analytic continuation with respect to $r$ from the positive real axis to the strip region 
$\{ r\in \mathbf{C} \, | \, -a<\mathrm{Im} (r) < a\}$.
Hence, the function $\mathcal{F}[\psi] (\sqrt{\lambda }\omega )$ of $\lambda $ has an analytic continuation from the 
positive real axis to the Riemann surface defined by
$P(a) = \{ \lambda \, | \, -a < \mathrm{Im} (\sqrt{\lambda }) < a \}$ with a branch point at the origin, see Fig.\ref{fig2}.
The region $P_1(a) = \{ \lambda \, | \, -a < \mathrm{Im} (\sqrt{\lambda }) < 0 \}$ is referred to as the first Riemann sheet,
and $P_2(a) = \{ \lambda \, | \, 0 < \mathrm{Im} (\sqrt{\lambda }) < a \}$ is referred to as the second Riemann sheet.
$P_1(a) \subset P_1 = \{ \lambda \, | \, -2\pi < \mathrm{arg} (\lambda ) < 0\}$ and
$P_2(a) \subset P_2 = \{ \lambda \, | \, 0 < \mathrm{arg} (\lambda ) < 2\pi\}$. 
They are connected with each other at the positive real axis $(\mathrm{arg} (\lambda ) = 0)$.
Therefore, the resolvent $(\lambda -H)^{-1}$ defined on the first Riemann sheet has an analytic continuation
to the second Riemann sheet through the positive real axis as
\begin{equation}
A(\lambda )i (\psi) (x) = 
\left\{ \begin{array}{ll}
(\lambda -H)^{-1}\psi (x) & (\lambda \in P_1), \\[0.2cm]
\displaystyle (\lambda -H)^{-1}\psi (x)+ 
\frac{\pi \sqrt{-1}}{(2\pi)^{m/2}}\sqrt{\lambda }^{m-2}\!\!\int_{S^{m-1}}\! \!
 e^{\sqrt{-1}\sqrt{\lambda }x \cdot \omega } \mathcal{F}[\psi] (\sqrt{\lambda }\omega ) d\omega
 & (\lambda \in P_2(a) ).\\
\end{array} \right.
\label{3-12}
\end{equation}
In particular, the space $X(\Omega ) = L^2(\mathbf{R}^m, e^{2a|x|}dx)$ satisfies (X1) to (X4) with $I=(0,\infty)$ and 
$\Omega = P_2(a)$.
To verify (X5), put $\sqrt{\lambda }=b_1 + \sqrt{-1}b_2$ with $-a < b_2 < a$.
Then, we obtain
\begin{eqnarray*}
|\mathcal{F}[\psi](\sqrt{\lambda }\omega )|^2
&\leq & \frac{1}{(2\pi)^{m}} \left| \int_{\mathbf{R}^m}\! \psi (x) e^{-\sqrt{-1}x \cdot \sqrt{\lambda }\omega } dx \right|^2 \\
&\leq & \frac{1}{(2\pi)^{m}} \left| \int_{\mathbf{R}^m}\! |\psi (x)| e^{|b_2x|} dx \right|^2 \\
&\leq & \frac{1}{(2\pi)^{m}} \left| \int_{\mathbf{R}^m}\! |\psi (x)| e^{(|b_2| - 2a)|x|}e^{2a|x|} dx \right|^2 \\
&\leq & \frac{1}{(2\pi)^{m}} \int_{\mathbf{R}^m}\!|\psi (x)|^2 e^{2a|x|} dx \cdot
           \int_{\mathbf{R}^m}\! e^{2(|b_2| - 2a)|x|}e^{2a|x|} dx \\
&=& \frac{1}{(2\pi)^{m}} || \psi ||^2_{X(\Omega )} \cdot \int_{\mathbf{R}^m}\! e^{2(|b_2| - a)|x|}dx.
\end{eqnarray*}
This proves that $\mathcal{F}[\psi] (\sqrt{\lambda }\omega )$ tends to zero 
uniformly in $\omega \in S^{m-1}$ as $\psi \to 0$ in $X(\Omega )$.
By using this fact and Eq.(\ref{3-12}), we can verify the assumption (X5).
Next, (X6) and (X7) are fulfilled with, for example, $Y = C^\infty_0 (\mathbf{R}^m)$,
the set of $C^\infty$ functions with compact support.
(X8) will be verified in the proof of Thm.3.1 below.
Since all assumptions are verified, the generalized spectrum $\mathcal{R}_\lambda $ of $H+K$ 
is well defined and Theorems 2.4, 2.6 and 2.7 hold (with $Z(\Omega ) = X(\Omega )$).
Let us show that (Z6) (for $Z(\Omega ) = X(\Omega )$) is satisfied and Thm.2.8 holds.
\\[0.2cm]
\textbf{Theorem \thedef.}
For any  $m=1,2,\cdots $, $i^{-1}K^\times A(\lambda )i$ is a compact operator on $L^2(\mathbf{R}^m, e^{2a|x|}dx)$.
In particular, the generalized spectrum $\hat{\sigma }(T)$ on the Riemann surface $P(a)$ 
consists only of a countable number of generalized eigenvalues having finite multiplicities.
\\[0.2cm]
\textbf{Proof.}
The multiplication by $e^{-a|x|}$ is a unitary operator from $L^2 (\mathbf{R}^m)$ onto $L^2(\mathbf{R}^m, e^{2a|x|}dx)$.
Thus, it is sufficient to prove that $e^{a|x|}i^{-1}K^\times A(\lambda )ie^{-a|x|}$ is a compact operator on
$L^2(\mathbf{R}^m)$.

When $\lambda $ lies on the first Riemann sheet, this operator is given as
\begin{equation}
e^{a|x|}i^{-1}K^\times A(\lambda )ie^{-a|x|}\psi (x) = e^{a|x|}V(x) (\lambda -H)^{-1}(e^{-a|x|}\psi (x)).
\label{3-15}
\end{equation}
The compactness of this type of operators is well known.
%
Indeed, it satisfies the Stummel condition for the compactness because of (\ref{4-11}).
See \cite{Rej, Wei} for the details.

Next, suppose that $\lambda \in I = (0,\infty)$; that is, $\lambda $ lies on the branch cut ($\mathrm{arg} (\lambda ) = 0$).
Let $L(X, X)$ be the set of continuous linear mappings on $X(\Omega )$ equipped with the usual norm topology.
A point $\lambda -\sqrt{-1}\varepsilon $ lies on the first Riemann sheet for small $\varepsilon >0$,
so that $i^{-1}K^\times A(\lambda -\sqrt{-1}\varepsilon ) i$ is compact.
It is sufficient to show that the sequence $\{ i^{-1}K^\times A(\lambda -\sqrt{-1}\varepsilon ) i\}_{\varepsilon >0}$
of compact operators converges to $i^{-1}K^\times A(\lambda) i$ as $\varepsilon \to 0$ in $L(X, X)$.
This immediately follows from the fact that $\{ i^{-1}K^\times A(\lambda) i\}_\lambda $ is a holomorphic family 
of operators; that is, $i^{-1}K^\times A(\lambda) i$ is holomorphic in $\lambda $ with respect to the topology of $L(X, X)$.
See Thm.3.12 of \cite{Chi} for the proof.

Finally, suppose that $\lambda $ lies on the second Riemann sheet.
In this case,
\begin{eqnarray*}
e^{a|x|}i^{-1}K^\times A(\lambda )ie^{-a|x|}\psi (x)
&=& e^{a|x|}V(x) (\lambda -H)^{-1}(e^{-a|x|}\psi (x)) \\
& & + \frac{\pi \sqrt{-1}}{(2\pi)^{m/2}}\sqrt{\lambda }^{m-2}\!\! \int_{S^{m-1}}\! \!
 e^{a|x|}V(x) e^{\sqrt{-1}x \cdot \sqrt{\lambda }\omega } \mathcal{F}[e^{-a|\,\cdot \,|}\psi] (\sqrt{\lambda }\omega ) d\omega .
\end{eqnarray*}
Since the first term in the right hand side above is compact, it is sufficient to prove that the mapping
\begin{equation}
\psi (x)\mapsto (K_2\psi)(x) := \int_{S^{m-1}}\! \!
 e^{a|x|}V(x) e^{\sqrt{-1}x \cdot \sqrt{\lambda }\omega } \mathcal{F}[e^{-a|\,\cdot \,|}\psi] (\sqrt{\lambda }\omega ) d\omega 
\end{equation}
on $L^2 (\mathbf{R}^m)$ is Hilbert-Schmidt.
This is rewritten as
\begin{equation}
(K_2\psi)(x) = \frac{1}{(2\pi)^{m/2}}\int_{\mathbf{R}^m}\! \int_{S^{m-1}}\! \!
 e^{a|x|}V(x) e^{\sqrt{-1}x \cdot \sqrt{\lambda }\omega } e^{-a|y|}\psi (y)e^{-\sqrt{-1} y \cdot \sqrt{\lambda }\omega } d\omega dy.
\end{equation}
This implies that $K_2$ is an integral operator with the kernel
\begin{equation}
k_2(x,y) :=  \frac{1}{(2\pi)^{m/2}} \int_{S^{m-1}}\! \!
 e^{a|x|}V(x) e^{-a|y|} e^{\sqrt{-1} (x-y) \cdot \sqrt{\lambda }\omega }d\omega .
\end{equation}
This is estimated as
\begin{eqnarray*}
|k_2(x,y)|^2 \leq \frac{1}{(2\pi)^{m}} e^{2a|x|} |V(x)|^2  e^{-2a|y|} \left| \int_{S^{m-1}}\!
e^{\sqrt{-1} (x-y) \cdot \sqrt{\lambda }\omega } d\omega   \right|^2.
\end{eqnarray*}
Putting $\sqrt{\lambda } = b_1 + \sqrt{-1}b_2$ with $0<b_2<a$ yields
\begin{eqnarray*}
|k_2(x,y)|^2
&\leq & \frac{1}{(2\pi)^{m}} e^{2a|x|} |V(x)|^2  e^{-2a|y|} \left( \int_{S^{m-1}}\! e^{b_2|x-y| } d\omega  \right)^2 \\
&\leq & \frac{1}{(2\pi)^{m}} e^{2a|x|} |V(x)|^2  e^{-2a|y|} e^{2b_2|x| + 2b_2 |y|} \mathrm{vol}(S^{m-1})^2 \\
&\leq & \frac{1}{(2\pi)^{m}} e^{4a|x|} |V(x)|^2  e^{2(b_2 - a)|x|}e^{2(b_2 - a)|y|} \mathrm{vol}(S^{m-1})^2 .
\end{eqnarray*}
Since $b_2 - a < 0$ and $e^{2a|x|} V(x)\in L^2(\mathbf{R}^m)$, we obtain $\int\!\int\! |k_2(x,y)|^2dxdy < \infty $,
which proves that $K_2$ is a Hilbert-Schmidt operator on $L^2(\mathbf{R}^m)$. \hfill $\blacksquare$


\subsection{Dilation analytic potentials}

In this subsection, we consider a dilation analytic potential.
At first, we define the van Winter space.

By using the polar coordinates, we denote a point $x\in \mathbf{R}^m$ as $x = r\omega ,\,\, r>0, \omega \in S^{m-1}$,
where $S^{m-1}$ is the $m-1$ dimensional unit sphere.
A function defined on $\mathbf{R}^m$ is denoted by $f(x) = f(r, \omega )$.
Any $f\in L^2(\mathbf{R}^m)$ satisfies
\begin{equation}
\int_{\mathbf{R}^m}\! |f(x)|^2 dx = \int_{S^{m-1}}\! d\omega \int^\infty_{0}\! |f(r, \omega )|^2 r^{m-1} dr < \infty.   
\end{equation}
Let $-\pi/2<\beta < \alpha < \pi /2$ be fixed numbers and $\mathcal{S}_{\beta, \alpha }  = \{ z\in \mathbf{C} \, | \, 
\beta < \mathrm{arg} (z) < \alpha \}$ an open sector.
Let $G (\beta, \alpha )$ be a vector space over $\mathbf{C}$ consisting of complex-valued functions
$f(z, \omega )$ on $\mathcal{S}_{\beta, \alpha }  \times S^{m-1}$ satisfying following conditions:
\\[0.2cm]
(i) $f(z, \omega )$ is holomorphic in $z\in \mathcal{S}_{\beta, \alpha } $.
\\
(ii) Put $z = re^{\sqrt{-1}\theta }$ with $r>0,\, \theta \in \mathbf{R}$. Then,
\begin{equation}
\sup_{\beta <\theta <\alpha } \int_{S^{m-1}}\! d\omega 
        \int^\infty_{0}\! |f(re^{\sqrt{-1}\theta }, \omega )|^2 r^{m-1} dr < \infty.
\end{equation}
It is known that when $f\in G(\beta, \alpha )$, boundary values $f(re^{\sqrt{-1}\alpha  }, \omega )$
and $f(re^{\sqrt{-1}\beta }, \omega )$ exist in $L^2$-sense.
By the inner product defined to be 
\begin{equation}
(f, g)_{\beta, \alpha } = \int_{S^{m-1}}\! d\omega \int^{\infty}_{0}\! \left( f(re^{\sqrt{-1}\alpha  }, \omega )
\overline{g(re^{\sqrt{-1}\alpha }, \omega )} + f(re^{\sqrt{-1}\beta  }, \omega )
\overline{g(re^{\sqrt{-1}\beta }, \omega )} \right) r^{m-1}dr,  
\end{equation}
$G(\beta, \alpha )$ becomes a Hilbert space.
In particular, when $\beta \leq 0 \leq \alpha $, $G(\beta, \alpha )$ is a dense subspace of $L^2(\mathbf{R}^m)$ and the topology of
$G(\beta, \alpha )$ is stronger than that of $L^2(\mathbf{R}^m)$.
The space $G(\beta, \alpha )$ was introduced by van Winter \cite{Win, Win2}.
In his papers, it is proved that the Fourier transform is a unitary mapping from $G(\beta, \alpha )$ onto $G(-\alpha , -\beta)$.
For a function $f(r, \omega ) \in G(\beta , \alpha )$ with $\beta \leq 0 \leq \alpha $, define the function $f^*$ by
$f^*(r, \omega ) = \overline{f(\overline{r}, \omega )}$.
Then, $f^* \in G(-\alpha , -\beta)$ and $f \mapsto f^*$ is a unitary mapping.
Let $\hat{f}$ be the Fourier transform of $f$.
It turns out that $(\hat{f})^*$ is in $G(\beta, \alpha )$ which is expressed as
\begin{equation}
(\hat{f})^*(re^{\sqrt{-1}\theta }, \omega )
 = \frac{e^{\sqrt{-1}m\theta }}{(2\pi)^{m/2}} \int_{S^{m-1}}\! d\omega ' \int^{\infty}_{0}\!
e^{\sqrt{-1}ry \omega \cdot \omega '} \overline{f(e^{\sqrt{-1}\theta }y, \omega ')}y^{m-1}dy.  
\end{equation}

Fix positive numbers $0<\beta < \alpha < \pi/2$.
We consider two spaces $G(-\alpha , \alpha )$ and $G(-\alpha, -\beta)$.
Since $G(-\alpha , \alpha )$ is dense in $L^2 (\mathbf{R}^m)$ and $G(-\alpha, -\beta )$,
we obtain the following diagram, in which the embedding $j$ will be defined later.

\begin{eqnarray*}
\begin{array}{ccccc}
G(-\alpha  , \alpha ) & \subset & L^2 (\mathbf{R}^m) & \subset & G(-\alpha  , \alpha )'\\
\rotatebox{90}{$\supset$} & & & & \rotatebox{90}{$\subset$} \\
G(-\alpha, -\beta) &  &\longrightarrow  & & jG(-\alpha, -\beta )\\
\end{array}
\end{eqnarray*}

In what follows, suppose that a potential $V : \mathbf{R}^m \to \mathbf{C}$ is an element of $G(-\alpha ,\alpha )$.
In the literature, such a potential is called a dilation analytic potential.
\\[0.2cm]
\textbf{Theorem \thedef.} Suppose $V\in G(-\alpha ,\alpha )$, $m=1,2,3$ and $0<\beta < \alpha < \pi/2$.
For $H + K = -\Delta + V$ and $\mathcal{H} = L^2 (\mathbf{R}^m)$, 
put $\Omega   = \{ \lambda \, | \, 0< \mathrm{arg} (\lambda ) < 2\alpha \}$, $I=(0, \infty)$, 
$X(\Omega ) = G(-\alpha  , \alpha )$ and $Z(\Omega ) = G(-\alpha, -\beta )$.
Then, they satisfy the assumptions (X1) to (X7) and (Z1) to (Z6).
\\[0.2cm]
\textbf{Corollary \thedef.} For each $\phi \in X(\Omega )$, the generalized resolvent 
$\mathcal{R}_\lambda i \phi$ of $-\Delta + V$ is an $X(\Omega )'$-valued meromorphic function defined on 
\begin{eqnarray*}
\hat{\Omega } = \{ \lambda \, | \, -2\pi-2 \alpha < \mathrm{arg} (\lambda ) < 2\alpha \}.
\end{eqnarray*}
The generalized spectrum consists of a countable number of generalized eigenvalues having finite multiplicities.
\\[0.2cm]
\textbf{Proof.}
The resolvent $(\lambda -T)^{-1}$ in the usual sense is meromorphic on 
$-2 \pi < \mathrm{arg}(\lambda ) < 0$.
Thm.2.4, 2.8 and 3.2 show that it has a meromorphic continuation from the sector 
$-\varepsilon < \mathrm{arg}(\lambda ) < 0$ to $\Omega $.
A similar argument also proves that it has a meromorphic continuation from the sector 
$-2 \pi + \varepsilon < \mathrm{arg}(\lambda ) < -2\pi$ to $-2 \pi < \mathrm{arg}(\lambda ) < -2\pi - 2\alpha $.
\hfill $\blacksquare$
\\[0.2cm]
Let us prove Thm.3.2.
The assumptions (X1), (X2), (X3), (Z1) and (Z2) are trivial.
\\[0.2cm]
\textbf{Proof of (X4), (X5).}
Let us calculate the analytic continuation of the resolvent $(\lambda -H)^{-1}\phi$.
For $\phi, \psi\in G(-\alpha ,\alpha )$, we have
\begin{eqnarray*}
((\lambda -H)^{-1}\phi, \psi) 
&=& \int_{\mathbf{R}^m}\! \frac{1}{\lambda -|\xi|^2}\mathcal{F}[\phi](\xi)\overline{\mathcal{F}[\psi](\xi)} d\xi \\
&=& \int_{S^{m-1}}\! d\omega \int^{\infty}_{0}\! \frac{1}{\lambda -y^2} \hat{\phi}(y, \omega )(\hat{\psi})^* (y, \omega )y^{m-1}dy,
\end{eqnarray*}
which is holomorphic on the first Riemann sheet $-2\pi < \mathrm{arg}(\lambda ) < 0$.
Since $\hat{\phi}, (\hat{\psi})^*\in G(-\alpha ,\alpha )$, Cauchy theorem yields
\begin{eqnarray*}
((\lambda -H)^{-1}\phi, \psi) =\int_{S^{m-1}}\! d\omega \int^{\infty}_{0}\! 
\frac{1}{\lambda -y^2e^{2\sqrt{-1}\alpha  }} \hat{\phi}(ye^{\sqrt{-1}\alpha  }, \omega )
        (\hat{\psi})^* (ye^{\sqrt{-1}\alpha  }, \omega )e^{\sqrt{-1}m\alpha  }y^{m-1}dy.
\end{eqnarray*}
This implies that $((\lambda -H)^{-1}\phi, \psi)$ is holomorphic on the sector $0\leq \mathrm{arg} (\lambda )< 2\alpha$.
When $\phi \to 0$ in $G(-\alpha , \alpha )$, then $\hat{\phi} \to 0$ in $G(-\alpha , \alpha )$.
Thus the above quantity tends to zero.
Similarly, $\psi \to 0$ in $G(-\alpha , \alpha )$ implies that $((\lambda -H)^{-1}\phi, \psi)$ tends to zero.
Therefore, the analytic continuation $A(\lambda ) : iX(\Omega ) \to X(\Omega )'$ of $(\lambda -H)^{-1}$
is given by
\begin{equation}
\langle A(\lambda )i\phi \,|\,\psi \rangle =\int_{S^{m-1}}\! d\omega \int^{\infty}_{0}\! 
\frac{1}{\lambda -y^2e^{2\sqrt{-1}\alpha  }} \hat{\phi}(ye^{\sqrt{-1}\alpha  }, \omega )
        (\hat{\psi})^* (ye^{\sqrt{-1}\alpha  }, \omega )e^{\sqrt{-1}m\alpha  }y^{m-1}dy,
\label{3-19}
\end{equation}
which is separately continuous in $\phi$ and $\psi \in G(-\alpha ,\alpha )$.
This confirms the assumptions (X4) and (X5) with $\Omega = \{ \lambda \, | \, 0< \mathrm{arg}(\lambda ) < 2\alpha \}$
and $I = (0, \infty)$.
\\[0.2cm]
\textbf{Proof of (X6), (X7).}
Let $\hat{Y}$ be the set of functions $\hat{f}\in G(-\alpha ,\alpha )$ such that $|\xi|^{2}\hat{f}(\xi) \in G(-\alpha ,\alpha )$.
Then, the inverse Fourier transform $Y$ of $\hat{Y}$ satisfies the assumption (X6).
Since $V\in G(-\alpha , \alpha ) \subset L^2(\mathbf{R}^m)$, Kato theorem proves that $V$ is $H$-bounded when $m=1,2,3$.
A proof of $V(Y) \subset G(-\alpha ,\alpha )$ is straightforward.
\\[0.2cm]
\textbf{Proof of (Z3).}
By the definition, the canonical inclusion $i : \mathcal{H} \to X(\Omega )'$ is defined through
\begin{eqnarray*}
\langle i\phi \,|\, \psi \rangle = (\phi, \psi)
 = \int_{S^{m-1}}\! d\omega \int^{\infty}_{0}\! \phi (r, \omega )\overline{\psi (r, \omega )}r^{m-1}dr,
\end{eqnarray*}
for $\phi \in \mathcal{H}$ and $\psi \in X(\Omega ) = G(-\alpha ,\alpha )$.
When $\phi \in X(\Omega )$, this is rewritten as
\begin{eqnarray*}
\langle i\phi \,|\, \psi \rangle
 = \int_{S^{m-1}}\! d\omega \int^{\infty}_{0}\! 
\phi (re^{\sqrt{-1}\theta }, \omega )\psi^* (re^{\sqrt{-1}\theta }, \omega )e^{\sqrt{-1}m\theta }r^{m-1}dr,
\end{eqnarray*}
for any $-\alpha \leq \theta \leq \alpha $.
The right hand side exists even for $\phi \in G(-\alpha , -\beta)$ when $-\alpha  \leq \theta \leq -\beta$.
Hence, we define the embedding $j : Z(\Omega ) \to X(\Omega )'$ to be 
\begin{eqnarray}
\langle j\phi \,|\, \psi \rangle
 = \int_{S^{m-1}}\! d\omega \int^{\infty}_{0}\! 
\phi (re^{-\sqrt{-1}\alpha }, \omega )\psi^* (re^{-\sqrt{-1}\alpha }, \omega )e^{-\sqrt{-1}m\alpha }r^{m-1}dr,
\end{eqnarray}
for $\phi \in Z(\Omega ) = G(-\alpha , -\beta)$.
This gives a continuous extension of $i : X(\Omega ) \to X(\Omega )'$.
In what follows, $j$ is denoted by $i$ for simplicity.
\\[0.2cm]
\textbf{Proof of (Z4), (Z5).}
The right hand side of Eq.(\ref{3-19}) is well-defined for $\phi \in Z(\Omega )$ and $\lambda \in \Omega $
because $\hat{\phi} \in G(\beta, \alpha )$ if $\phi \in G(-\alpha ,-\beta)$.
This gives an extension of  $A(\lambda ) : iX(\Omega ) \to X(\Omega )'$ to $A(\lambda ) : iZ(\Omega ) \to X(\Omega )'$.
(Z5) is verified as follows: for $\phi \in G(-\alpha , -\beta )$ and $-\alpha <\theta <-\beta$,
\begin{eqnarray}
& & |(VA(\lambda )i \phi ) (re^{\sqrt{-1}\theta }, \omega )|^2 \nonumber \\
&= & \frac{1}{(2\pi)^m} \left| \int\! d\omega' \int\! \frac{V(re^{\sqrt{-1}\theta }, \omega )}{\lambda -y^2e^{-2\sqrt{-1}\theta}}
       e^{\sqrt{-1}ry \omega \cdot \omega '} \hat{\phi}(ye^{-\sqrt{-1}\theta }, \omega ')y^{m-1}dy \right|^2 \nonumber\\
&\leq & \frac{1}{(2\pi)^m} \int\! d\omega' \int\! \frac{|V(re^{\sqrt{-1}\theta },\omega )|^2}{|\lambda -y^2e^{-2\sqrt{-1}\theta} |^2}
       y^{m-1}dy \cdot \int\! d\omega' \int\! |\hat{\phi}(ye^{-\sqrt{-1}\theta }, \omega ')|^2 y^{m-1}dy\nonumber \\
&\leq &C_1  \int\! d\omega' \int\! \frac{|V(re^{\sqrt{-1}\theta },\omega )|^2}{|\lambda -y^2e^{-2\sqrt{-1}\theta }|^2} y^{m-1}dy,
\label{3-34}
\end{eqnarray} 
where 
\begin{eqnarray*}
C_1 = \frac{1}{(2\pi)^m} \sup_{-\alpha \leq \theta \leq -\beta} 
\int\! d\omega' \int\! |\hat{\phi}(ye^{-\sqrt{-1}\theta }, \omega ')|^2 y^{m-1}dy,
\end{eqnarray*}
which exists because $\hat{\phi} \in G(\beta, \alpha )$.
Hence, we obtain
\begin{eqnarray*}
& & \int\!d\omega \int\! |VA(\lambda )i \phi (re^{\sqrt{-1}\theta }, \omega )|^2 r^{m-1}dr \\
& \leq & C_1  \int\!d\omega \int\! |V(re^{\sqrt{-1}\theta },\omega )|^2  r^{m-1}dr \cdot
          \int\! d\omega' \int\! \frac{1}{|\lambda -y^2e^{-2\sqrt{-1}\theta }|^2} y^{m-1}dy.
\end{eqnarray*}
Since $V\in G(-\alpha , \alpha )$ and $m=1,2,3$, this has an upper bound which is independent of 
$-\alpha  < \theta < -\beta$.
This shows $VA(\lambda )i \phi \in G(-\alpha , -\beta)$.
Let $|| \cdot ||_{\beta, \alpha } $ be the norm on $G(\beta, \alpha )$.
Since $C_1 \leq || \hat{\phi} ||_{\beta, \alpha }/(2\pi)^{m} = || \phi ||_{-\alpha , -\beta}/(2\pi)^{m}$,
it immediately follows that $V A(\lambda )i : G(-\alpha , -\beta) \to G(-\alpha , -\beta)$ is continuous.
\\[0.2cm]
\textbf{Proof of (Z6).}
Let $P :  G(-\alpha , -\beta) \to G(-\alpha , -\beta)$ be a continuous linear operator.
For fixed $-\alpha \leq \theta \leq -\beta$, $f(re^{\sqrt{-1}\theta }, \omega )$ and 
$(Pf)(re^{\sqrt{-1}\theta }, \omega )$ are regarded as functions of $L^2(\mathbf{R}^m)$.
Thus the mapping $f(re^{\sqrt{-1}\theta }, \omega ) \mapsto (Pf)(re^{\sqrt{-1}\theta }, \omega ) $
defines a continuous linear operator on $L^2 (\mathbf{R}^m)$, which is denoted by $P_\theta $
(this was introduced by van Winter \cite{Win}).
To verify (Z6), we need the next lemma.
\\[0.2cm]
\textbf{Lemma \thedef.}
If $P_{-\alpha }$ and $P_{-\beta}$ are compact operators on $L^2(\mathbf{R}^m)$,
then $P$ is a compact operator on $G(-\alpha , -\beta)$.
\\[0.2cm]
\textbf{Proof.} 
Let $B\subset G(-\alpha , -\beta)$ be a bounded set.
Since the topology of $L^2(\mathbf{R}^m)$ is weaker than that of $G(-\alpha , -\beta)$,
$\{ f(re^{-\sqrt{-1}\alpha }, \omega )\}_{f\in B}$ and 
$\{ f(re^{-\sqrt{-1}\beta }, \omega )\}_{f\in B}$ are bounded sets of $L^2(\mathbf{R}^m)$.
Since $P_{-\alpha }$ is compact, there exists a sequence 
$\{ g_j\}^\infty_{j=1} \subset B$ such that $(P_{-\alpha }g_j)(re^{-\sqrt{-1}\alpha }, \omega )$
converges in $L^2(\mathbf{R}^m)$.
Since $P_{-\beta }$ is compact, there exists a sequence 
$\{ h_j\}^\infty_{j=1} \subset \{ g_j\}^\infty_{j=1}$ such that $(P_{-\beta }h_j)(re^{-\sqrt{-1}\beta }, \omega )$
converges in $L^2(\mathbf{R}^m)$.
By the definition of the norm of $G(-\alpha , -\beta)$, $Ph_j$ converges in $G(-\alpha , -\beta)$. \hfill $\blacksquare$

In order to verify (Z6), it is sufficient to show that the mapping $\phi (re^{\sqrt{-1}\theta }, \omega ) \mapsto$
 $(VA(\lambda )i\phi)(re^{\sqrt{-1}\theta }, \omega )$ on $L^2(\mathbf{R}^m)$ is compact for $\theta =-\alpha ,-\beta$.
The $(VA(\lambda )i\phi)(re^{-\sqrt{-1}\alpha }, \omega )$ is given by
\begin{eqnarray*}
& & (VA(\lambda )i\phi)(re^{-\sqrt{-1}\alpha }, \omega ) \\
&=&  \frac{e^{\sqrt{-1}m\alpha }}{(2\pi)^{m/2}} \!\int\!\! d\omega_1 \!\!\int\!
\frac{V(re^{-\sqrt{-1}\alpha }, \omega )}{\lambda -y_1^2e^{2\sqrt{-1}\alpha }} e^{\sqrt{-1}ry_1 \omega \cdot \omega_1}
       \hat{\phi}(y_1e^{\sqrt{-1}\alpha }, \omega_1)y_1^{m-1}dy_1 \\
&=& \!\!\!\! \frac{1}{(2\pi)^{m}} \!\iint \!\! d\omega_1d\omega _2 \!\!\iint \!
\frac{V(re^{-\sqrt{-1}\alpha }, \omega )}{\lambda -y_1^2e^{2\sqrt{-1}\alpha }} e^{\sqrt{-1}ry_1 \omega \cdot \omega_1}
 e^{-\sqrt{-1} y_1y_2 \omega _1 \cdot \omega _2}
   \phi(y_2e^{-\sqrt{-1}\alpha }, \omega_2)y_1^{m-1}y_2^{m-1}dy_1dy_2. 
\end{eqnarray*}
This is an integral operator with the kernel
\begin{eqnarray*}
K(r,\omega ;y_2,\omega _2)
=\frac{1}{(2\pi)^{m}} \!\int \!\!d\omega _1 \!\!\int \!
\frac{V(re^{-\sqrt{-1}\alpha }, \omega )}{\lambda -y_1^2e^{2\sqrt{-1}\alpha }} e^{-\sqrt{-1}y_1 \omega_1 (y_2\omega _2-r\omega )}
y_1^{m-1}dy_1. 
\end{eqnarray*}
Put $f(y_1, \omega _1) = (\lambda -y_1^2e^{2\sqrt{-1}\alpha })^{-1}$, which is in $L^2(\mathbf{R}^m)$ when $m=1,2,3$.
Then, we obtain
\begin{eqnarray*}
|K(r,\omega ;y_2,\omega _2)|^2 \leq \frac{1}{(2\pi)^m} |V(re^{-\sqrt{-1}\alpha }, \omega )|^2 \cdot 
 |\hat{f}(y_2\omega _2 - r\omega )|^2.
\end{eqnarray*}
Putting $x=r\omega , \xi = y_2\omega _2$ yields
\begin{eqnarray*}
\iint |K(r,\omega ;y_2,\omega _2)|^2 dxd\xi \leq \frac{1}{(2\pi)^m} \int |V(xe^{-\sqrt{-1}\alpha })|^2 dx \cdot 
\int |\hat{f}(\xi -x )|^2 d\xi.
\end{eqnarray*}
This quantity exists because $\hat{f}\in L^2(\mathbf{R}^m)$ and $V\in G(-\alpha , \alpha )$.
Therefore, the mapping $\phi (re^{\sqrt{-1}\theta }, \omega ) \mapsto (VA(\lambda )i\phi)(re^{\sqrt{-1}\theta }, \omega )$ 
is a Hilbert-Schmidt operator for $\theta =-\alpha $.
The other case $\theta =-\beta$ is proved in the same way.
Now the proof of Thm.3.2 is completed. \hfill $\blacksquare$


\section{One dimensional Schr\"{o}dinger operators}

In this section, one dimensional Schr\"{o}dinger operators (ordinary differential operators) are considered.
A holomorphic function $\mathbb{D}(\lambda ,\nu_1, \nu_2)$, zeros of which coincide with generalized eigenvalues,
is constructed.
It is proved that the function $\mathbb{D}(\lambda ,\nu_1, \nu_2)$ is equivalent to the Evans function.


\subsection{Evans functions for exponentially decaying potentials}

Let us consider the one dimensional Schr\"{o}dinger operator on $L^2(\mathbf{R})$
\begin{equation}
-\frac{d^2}{dx^2} + V(x), \quad x\in \mathbf{R},
\label{4-1}
\end{equation}
with an exponentially decaying potential $e^{2a|x|}V(x) \in L^2(\mathbf{R})$ for some $a>0$.
A point $\lambda $ is a generalized eigenvalue if and only if $(id -  A(\lambda )K^\times)\mu = 0$ has a
nonzero solution $\mu$ in $X(\Omega )' = L^2(\mathbf{R}, e^{-2a|x|}dx)$.
For a one dimensional case, this equation is reduced to the integral equation
\begin{equation}
\mu(x, \lambda ) 
= -\frac{1}{2\sqrt{-\lambda }} \int^{\infty}_{x}\! e^{-\sqrt{-\lambda } (y-x)}V(y)\mu(y, \lambda )dy 
-\frac{1}{2\sqrt{-\lambda }} \int^{x}_{-\infty}\! e^{\sqrt{-\lambda } (y-x)}V(y)\mu(y, \lambda )dy,
\label{3-27}
\end{equation}
where $\lambda $ lies on the Riemann surface $P(a) = \{ \lambda \, | \, -a < \mathrm{Im}(\sqrt{\lambda }) < a\}$.
It is convenient to rewrite it as a differential equation.
Due to Thm.2.6, $\mu$ satisfies $(\lambda -T^\times) \mu = 0$.
This is a differential equation on $ L^2(\mathbf{R}, e^{-2a|x|}dx)$ of the form 
\begin{equation}
\left( \frac{d^2}{dx^2} + \lambda - V(x) \right) \mu(x, \lambda ) = 0,\quad 
\mu (\, \cdot \, , \lambda ) \in L^2(\mathbf{R}, e^{-2a|x|}dx).
\label{4-3}
\end{equation}
Since $V$ decays exponentially, for any solutions $\mu$, there exists a constant $C$ such that 
$|\mu (x, \lambda )| \leq C e^{|\mathrm{Im} \sqrt{\lambda }| \cdot |x|}$ (see \cite{Cod}, Chap.3).
In particular, any solutions satisfy $\mu \in L^2(\mathbf{R}, e^{-2a|x|}dx)$ when $\lambda \in P(a)$.
At first, we will solve the differential equation (\ref{4-3}) for any $\lambda $,
which gives a candidate $\mu (x, \lambda )$ of a generalized eigenfunction.
Substituting the candidate into Eq.(\ref{3-27}) determines a generalized eigenvalue $\lambda $ 
and a true generalized eigenfunction $\mu$ associated with $\lambda $.
For this purpose, we define the function $D(x, \lambda , \mu)$ to be
\begin{eqnarray*}
& & D(x, \lambda , \mu) = \mu - A(\lambda )K^\times \mu =  \nonumber \\ 
& & \mu(x, \lambda ) 
+ \frac{1}{2\sqrt{-\lambda }} \int^{\infty}_{x}\! e^{-\sqrt{-\lambda } (y-x)}V(y)\mu(y, \lambda )dy 
+\frac{1}{2\sqrt{-\lambda }} \int^{x}_{-\infty}\! e^{\sqrt{-\lambda } (y-x)}V(y)\mu(y, \lambda )dy.
\end{eqnarray*}
Let $\mu_1$ and $\mu_2$ be two linearly independent solutions of Eq.(\ref{4-3}).
Define the function $\mathbb{D}(\lambda , \mu_1, \mu_2)$ to be
\begin{equation}
\mathbb{D}(\lambda , \mu_1, \mu_2) = \mathrm{det} \left(
\begin{array}{@{\,}cc@{\,}}
D(x, \lambda , \mu_1) & D(x, \lambda , \mu_2) \\
D'(x, \lambda , \mu_1) & D'(x, \lambda , \mu_2) 
\end{array}
\right),
\label{4-4}
\end{equation}
where $D'(x, \lambda , \mu)$ denotes the derivative with respect to $x$.
It is easy to verify that $\mathbb{D}(\lambda , \mu_1, \mu_2)$ is independent of $x$.
Next, let $\nu_1$ and $\nu_2$ be two solutions of Eq.(\ref{4-3}) satisfying the initial conditions
\begin{equation}
\left(
\begin{array}{@{\,}c@{\,}}
\nu_1 (0,\lambda ) \\
\nu_1 '(0,\lambda )
\end{array}
\right) = \left(
\begin{array}{@{\,}c@{\,}}
1 \\
0
\end{array}
\right), \quad \left(
\begin{array}{@{\,}c@{\,}}
\nu_2 (0,\lambda ) \\
\nu_2 '(0,\lambda )
\end{array}
\right) = \left(
\begin{array}{@{\,}c@{\,}}
0 \\
1
\end{array}
\right),
\label{4-5}
\end{equation}
respectively, where an initial time is fixed arbitrarily.
Because of the linearity, we obtain
\begin{equation}
\mathbb{D}(\lambda , \mu_1, \mu_2) = \mathrm{det} \left(
\begin{array}{@{\,}cc@{\,}}
\mu_1(0, \lambda ) & \mu_2(0, \lambda ) \\
\mu_1'(0, \lambda ) & \mu_2'(0, \lambda )
\end{array}
\right) \cdot \mathbb{D}(\lambda , \nu_1, \nu_2) .
\label{4-6}
\end{equation}
Therefore, it is sufficient to investigate properties of $\mathbb{D}(\lambda , \nu_1, \nu_2)$.
The main theorem in this section is stated as follows.
\\[0.2cm]
\textbf{Theorem \thedef.} For the operator (\ref{4-1}) with $e^{2a|x|}V\in L^2(\mathbf{R})$,
\\
(i) $\mathbb{D}(\lambda , \nu_1, \nu_2)$ is holomorphic in $\lambda \in P(a)$.
\\
(ii) $\mathbb{D}(\lambda , \nu_1, \nu_2) = 0$ if and only if $\lambda $ is a generalized eigenvalue of (\ref{4-1}).
\\
(iii) When $\lambda $ lies on the first Riemann sheet ($-2\pi < \mathrm{arg} (\lambda ) < 0$),
\begin{equation}
\mathbb{D}(\lambda , \nu_1, \nu_2) = \frac{1}{2\sqrt{-\lambda }} \mathbb{E}(\lambda ),
\end{equation}
where $ \mathbb{E}(\lambda )$ is the Evans function defined in Sec.1.
\\[0.2cm]
\textbf{Corollary \thedef.}
The Evans function $\mathbb{E}(\lambda )$ has an analytic continuation from the plane
$-2\pi < \mathrm{arg} (\lambda ) < 0$ to the Riemann surface $P(a)$,
whose zeros give generalized eigenvalues.
\\[0.2cm]
\textbf{Proof.} (i) Since initial conditions (\ref{4-5}) are independent of $\lambda $,
$\nu_1(x, \lambda )$ and $\nu_2(x, \lambda )$ are holomorphic in $\lambda \in \mathbf{C}$.
Since $A(\lambda )$ is holomorphic in $\lambda \in P(a)$, $D(x, \lambda , \nu_i) = (id - A(\lambda )K^\times) \nu_i(x, \lambda )$
is holomorphic in $\lambda \in P(a)$ for $i=1,2$, which proves (i).

To prove (ii), suppose that $\lambda $ is a generalized eigenvalue.
Thus there exists $\mu \in L^2(\mathbf{R}, e^{-2a|x|}dx)$ such that $\mu - A(\lambda )K^\times \mu = 0$.
Since the generalized eigenfunction $\mu$ has to be a solution of the differential equation (\ref{4-3}),
there are numbers $C_1(\lambda )$ and $C_2(\lambda )$ such that $\mu = C_1(\lambda )\nu_1 + C_2(\lambda )\nu_2$.
Hence,
\begin{eqnarray*}
0& = & \mu - A(\lambda )K^\times \mu  = C_1(\lambda ) (\nu_1 - A(\lambda )K^\times \nu_1) +
C_2(\lambda ) (\nu_2 - A(\lambda )K^\times \nu_2) \\
&=& C_1(\lambda )D(x, \lambda , \nu_1) + C_2(\lambda )D(x, \lambda , \nu_2).
\end{eqnarray*}
This implies $\mathbb{D}(\lambda , \nu_1, \nu_2) = 0$.
Conversely, suppose that $\mathbb{D}(\lambda , \nu_1, \nu_2) = 0$.
Then, there are numbers $C_1(\lambda )$ and $C_2(\lambda )$ such that $(C_1 (\lambda ), C_2(\lambda )) \neq (0,0)$ and 
$C_1(\lambda )D(x, \lambda , \nu_1) + C_2(\lambda )D(x, \lambda , \nu_2) = 0$.
Putting $\mu =  C_1(\lambda )\nu_1 + C_2(\lambda )\nu_2$ provides $\mu - A(\lambda )K^\times \mu = 0$.
Due to the choice of the initial conditions of $\nu_1$ and $\nu_2$, we can show that $\mu (x, \lambda ) \nequiv 0$.
Hence, $\mu$ is a generalized eigenfunction associated with $\lambda $.

Finally, let us prove (iii).
Suppose that $-2\pi < \mathrm{arg} (\lambda ) < 0$, so that $\mathrm{Re}\sqrt{-\lambda } > 0$.
Let $\mu_+$ and $\mu_-$ be solutions of Eq.(\ref{4-3}) satisfying (\ref{1-1}), (\ref{1-2})
(for the existence of such solutions, see \cite{Cod}, Chap.3).
The Evans function is defined by $\mathbb{E}(\lambda ) = \mu_+(x, \lambda )\mu_-'(x, \lambda )
- \mu_+'(x, \lambda )\mu_-(x, \lambda )$.
Since $\mathbb{E}(\lambda )$ is independent of $x$, Eq.(\ref{4-6}) yields
\begin{equation}
\mathbb{D}(\lambda , \mu_+, \mu_-) = \mathrm{det} \left(
\begin{array}{@{\,}cc@{\,}}
\mu_+(0, \lambda ) & \mu_-(0, \lambda ) \\
\mu_+'(0, \lambda ) & \mu_-'(0, \lambda ) 
\end{array}
\right) \cdot \mathbb{D}(\lambda , \nu_1, \nu_2) = \mathbb{E}(\lambda ) \cdot \mathbb{D}(\lambda , \nu_1, \nu_2).
\end{equation}
Thus, it is sufficient to prove the equality
\begin{equation}
\mathbb{D}(\lambda , \mu_+, \mu_-) = \frac{1}{2\sqrt{-\lambda }}\mathbb{E}(\lambda )^2.
\label{4-9}
\end{equation}
At first, note that $\mu_{\pm}$ satisfy the integral equations
\begin{eqnarray*}
& & \mu_+ (x, \lambda ) = e^{-\sqrt{-\lambda }x} 
  - \frac{1}{2\sqrt{-\lambda }}\int^{\infty}_{x}\!\!e^{-\sqrt{-\lambda }(y-x)} V(y)\mu_+ (y, \lambda )dy
 +  \frac{1}{2\sqrt{-\lambda }} \int^{\infty}_{x}\!\!e^{\sqrt{-\lambda }(y-x)} V(y)\mu_+ (y, \lambda )dy, \\
& & \mu_- (x, \lambda ) = e^{\sqrt{-\lambda }x} 
+ \frac{1}{2\sqrt{-\lambda }} \int^{x}_{-\infty}\!\!e^{-\sqrt{-\lambda }(y-x)} V(y)\mu_- (y, \lambda )dy
 - \frac{1}{2\sqrt{-\lambda }} \int^{x}_{-\infty}\!\!e^{\sqrt{-\lambda }(y-x)} V(y)\mu_- (y, \lambda )dy.
\end{eqnarray*}
Substituting them into the definition of $D(x, \lambda , \mu)$ yields
\begin{eqnarray*}
\left\{ \begin{array}{l}
\displaystyle D(x, \lambda , \mu_+) =  e^{-\sqrt{-\lambda }x} 
    + \frac{1}{2\sqrt{-\lambda }}\int_{\mathbf{R}}\!e^{\sqrt{-\lambda }(y-x)} V(y)\mu_+ (y, \lambda )dy, \\
\displaystyle D'(x, \lambda , \mu_+) = -\sqrt{-\lambda } D(x, \lambda , \mu_+), \\
\displaystyle D(x, \lambda , \mu_-) =  e^{\sqrt{-\lambda }x} 
    + \frac{1}{2\sqrt{-\lambda }}\int_{\mathbf{R}}\!e^{-\sqrt{-\lambda }(y-x)} V(y)\mu_- (y, \lambda )dy, \\
\displaystyle D'(x, \lambda , \mu_-) = \sqrt{-\lambda } D(x, \lambda , \mu_+).
\end{array} \right.
\end{eqnarray*}
Next, since $\mathbb{E}(\lambda )$ is independent of $x$, we obtain
\begin{eqnarray*}
\mathbb{E}(\lambda ) &=& \lim_{x\to \infty}(\mu_+(x, \lambda )\mu_-'(x, \lambda ) - \mu_+'(x, \lambda )\mu_-(x, \lambda )) \\
&=& \lim_{x\to \infty} e^{-\sqrt{-\lambda }x} (\mu_-'(x, \lambda ) + \sqrt{-\lambda } \mu_-(x, \lambda )) \\
&=& 2 \sqrt{-\lambda } \left( 1+\frac{1}{2\sqrt{-\lambda }}\int_{\mathbf{R}}\!e^{-\sqrt{-\lambda }y} V(y)\mu_- (y, \lambda )dy\right),
\end{eqnarray*}
and 
\begin{eqnarray*}
\mathbb{E}(\lambda ) &=& \lim_{x\to -\infty}(\mu_+(x, \lambda )\mu_-'(x, \lambda ) - \mu_+'(x, \lambda )\mu_-(x, \lambda )) \\
&=& \lim_{x\to -\infty} e^{\sqrt{-\lambda }x} (\sqrt{-\lambda }\mu_+(x, \lambda ) -\mu_+'(x, \lambda )) \\
&=& 2 \sqrt{-\lambda } \left( 1+\frac{1}{2\sqrt{-\lambda }}\int_{\mathbf{R}}\!e^{\sqrt{-\lambda }y} V(y)\mu_+(y, \lambda )dy\right).
\end{eqnarray*}
They give
\begin{eqnarray*}
& & D(x, \lambda , \mu_+) =  \frac{e^{-\sqrt{-\lambda }x}}{2\sqrt{-\lambda }}\mathbb{E}(\lambda ),
\quad D(x, \lambda , \mu_-) =  \frac{e^{\sqrt{-\lambda }x}}{2\sqrt{-\lambda }}\mathbb{E}(\lambda ),\\
& & D'(x, \lambda , \mu_+) = -\frac{e^{-\sqrt{-\lambda }x} }{2}\mathbb{E}(\lambda ),
\quad D'(x, \lambda , \mu_-) = \frac{e^{\sqrt{-\lambda }x} }{2}\mathbb{E}(\lambda ).\\
\end{eqnarray*}
This proves Eq.(\ref{4-9}). \hfill $\blacksquare$


\subsection{Evans functions for dilation analytic potentials}

Let $\alpha $ and $\beta$ are positive numbers such that $0<\beta < \alpha < \pi/2$.
Let us consider the operator (\ref{4-1}) with a dilation analytic potential $V \in G(-\alpha ,\alpha )$.
For a one dimensional case, the equation $\mu = A(\lambda )K^\times \mu$ is written as
\begin{eqnarray*}
\mu (xe^{\sqrt{-1}\theta }, \lambda ) &=& 
-\frac{e^{\sqrt{-1}\theta }}{2\sqrt{-\lambda }} \int^{\infty}_{x}\! e^{-\sqrt{-\lambda }e^{\sqrt{-1}\theta }(y-x)}
  V(ye^{\sqrt{-1}\theta })\mu(ye^{\sqrt{-1}\theta }, \lambda )dy \\
& & -\frac{e^{\sqrt{-1}\theta }}{2\sqrt{-\lambda }} \int^{x}_{-\infty}\! e^{\sqrt{-\lambda }e^{\sqrt{-1}\theta } (y-x)}
V(ye^{\sqrt{-1}\theta })\mu(ye^{\sqrt{-1}\theta }, \lambda )dy,
\end{eqnarray*}
where $-\alpha < \theta < -\beta$.
Let $\mu$ be a solution of the differential equation (\ref{4-3}).
Since $V(x)$ is holomorphic in the sector $ \{ -\alpha < \mathrm{arg}(x) < \alpha \}$, 
so is the solution $\mu (x, \lambda )$.
Thus $\mu(xe^{\sqrt{-1}\theta }, \lambda )$ is well defined for $-\alpha < \theta < \alpha $.
For a solution $\mu$ of (\ref{4-3}), define the function $D(x, \lambda , \mu)$ to be 
\begin{eqnarray*}
D(x, \lambda , \mu) = \mu (xe^{\sqrt{-1}\theta }, \lambda )  
&+& \frac{e^{\sqrt{-1}\theta }}{2\sqrt{-\lambda }} \int^{\infty}_{x}\! e^{-\sqrt{-\lambda }e^{\sqrt{-1}\theta }(y-x)}
  V(ye^{\sqrt{-1}\theta })\mu(ye^{\sqrt{-1}\theta }, \lambda )dy \\
& & \quad +\frac{e^{\sqrt{-1}\theta }}{2\sqrt{-\lambda }} \int^{x}_{-\infty}\! e^{\sqrt{-\lambda }e^{\sqrt{-1}\theta } (y-x)}
V(ye^{\sqrt{-1}\theta })\mu(ye^{\sqrt{-1}\theta }, \lambda )dy.
\end{eqnarray*}
Then, the function $\mathbb{D}(\lambda , \mu_1, \mu_2) $ is defined by Eq.(\ref{4-4}).
By the same way as the previous subsection, it turns out that $\mathbb{D}(\lambda , \nu_1, \nu_2)=0 $
if and only if $\lambda $ is a generalized eigenvalue.

To define the Evans function, we further suppose that $V$ satisfies $\int_{\mathbf{R}}\! |V(xe^{\sqrt{-1}\theta })|dx < \infty$
for $-\alpha <\theta <-\beta$.
Then, there exist solutions $\mu_+$ and $\mu_-$ of Eq.(\ref{4-3}) satisfying
\begin{equation}
\left(
\begin{array}{@{\,}c@{\,}}
\mu_+(xe^{\sqrt{-1}\theta }, \lambda ) \\
\mu_+'(xe^{\sqrt{-1}\theta }, \lambda )
\end{array}
\right) e^{\sqrt{-\lambda }xe^{\sqrt{-1}\theta }} \to \left(
\begin{array}{@{\,}c@{\,}}
1 \\
-\sqrt{-\lambda }
\end{array}
\right), \quad (x\to \infty), 
\label{1-1b}
\end{equation}
and 
\begin{equation}
\left(
\begin{array}{@{\,}c@{\,}}
\mu_-(xe^{\sqrt{-1}\theta }, \lambda ) \\
\mu_-'(xe^{\sqrt{-1}\theta }, \lambda )
\end{array}
\right) e^{-\sqrt{-\lambda }xe^{\sqrt{-1}\theta }} \to \left(
\begin{array}{@{\,}c@{\,}}
1 \\
\sqrt{-\lambda }
\end{array}
\right), \quad (x\to -\infty), 
\label{1-2b}
\end{equation}
for some $-\alpha < \theta < -\beta$ (see \cite{Cod}, Chap.3).
For such $\theta $, define the Evans function $\mathbb{E}(\lambda )$ to be
\begin{equation}
\mathbb{E}(\lambda ) = e^{\sqrt{-1}\theta }\mu_+(xe^{\sqrt{-1}\theta }, \lambda ) \mu_-'(xe^{\sqrt{-1}\theta }, \lambda ) 
 - e^{\sqrt{-1}\theta }\mu_+'(xe^{\sqrt{-1}\theta }, \lambda ) \mu_-(xe^{\sqrt{-1}\theta }, \lambda ) ,
\end{equation}
which is independent of $x$ and $\theta $.
Then, we can prove the next theorem.
\\[0.2cm]
\textbf{Theorem \thedef.} Let $\hat{\Omega }$ be the region defined in Cor.3.3.
For the operator (\ref{4-1}) satisfying $V\in G(-\alpha , \alpha )$ and 
$\int_{\mathbf{R}}\! |V(xe^{\sqrt{-1}\theta })|dx < \infty$,
\\
(i) $\mathbb{D}(\lambda , \nu_1, \nu_2)$ is holomorphic in $\lambda \in \hat{\Omega }$.
\\
(ii) $\mathbb{D}(\lambda , \nu_1, \nu_2) = 0$ if and only if $\lambda $ is a generalized eigenvalue of (\ref{4-1}).
\\
(iii) When $\lambda $ lies on the first Riemann sheet ($-2\pi < \mathrm{arg} (\lambda ) < 0$),
\begin{equation}
\mathbb{D}(\lambda , \nu_1, \nu_2) = \frac{1}{2\sqrt{-\lambda }} \mathbb{E}(\lambda ).
\end{equation}
\textbf{Corollary \thedef.}
The Evans function $\mathbb{E}(\lambda )$ has an analytic continuation from the plane
$-2\pi < \mathrm{arg} (\lambda ) < 0$ to the Riemann surface $\hat{\Omega }$,
whose zeros give generalized eigenvalues.

The proofs are the same as before and omitted.
Note that (i) and (ii) hold without the assumption $\int_{\mathbf{R}}\! |V(xe^{\sqrt{-1}\theta })|dx < \infty$.


\subsection{Examples}

To demonstrate the theory developed so far, we give two solvable examples.
The first one is the simple one-dimensional potential well given by
\begin{equation}
V(x) = \left\{ \begin{array}{ll}
h & (a_1 < x < a_2), \\
0 & (\mathrm{otherwise}), \\
\end{array} \right.
\label{4-17}
\end{equation}
where $h, a_1, a_2 \in \mathbf{R}$.
Since $V(x)$ has compact support, it satisfies Eq.(\ref{4-11}) for any $a>0$.
Hence, we put
\begin{eqnarray*}
X(\Omega ) = \bigcap^\infty_{a\geq 0} L^2(\mathbf{R}, e^{2a|x|}dx) ,
\quad X(\Omega )' = \bigcup^\infty_{a\geq 0} L^2(\mathbf{R}, e^{-2a|x|}dx).
\end{eqnarray*}
They are equipped with the projective limit topology and the inductive limit topology, respectively.
Then, $X(\Omega )$ is a reflective Fr\'{e}chet space and the assumptions (X1) to (X8) are again 
verified in the same way as before.
In particular, $A(\lambda ), \mathcal{R}_\lambda $ and the generalized spectrum are defined on 
$\bigcup_{a\geq 1} P(a)$, which is the whole Riemann surface of $\sqrt{\lambda }$.
By applying Prop.2.9 with $Z_1(\Omega ) = L^2(\mathbf{R}^m, e^{2a|x|}dx)$ and 
$Z_2(\Omega ) = \bigcap^\infty_{a\geq 0} L^2(\mathbf{R}^m, e^{2a|x|}dx)$, Thm.3.1 immediately provides the next theorem,
which is true for any $m\geq 1$.
\\[0.2cm]
\textbf{Theorem \thedef.} Put $X(\Omega ) = \bigcap^\infty_{a\geq 0} L^2(\mathbf{R}^m, e^{2a|x|}dx)$.
Suppose that a potential $V$ satisfies Eq.(\ref{4-11}) for any $a\geq 0$.
Then, the generalized resolvent $\mathcal{R}_\lambda i \phi$ is a $X(\Omega )'$-valued meromorphic function
on the Riemann surface of $\sqrt{\lambda }$ for any $\phi \in X(\Omega )$.
In particular, the generalized spectrum consists of a countable number of generalized eigenvalues having finite multiplicities.

\begin{figure}
\begin{center}
\includegraphics{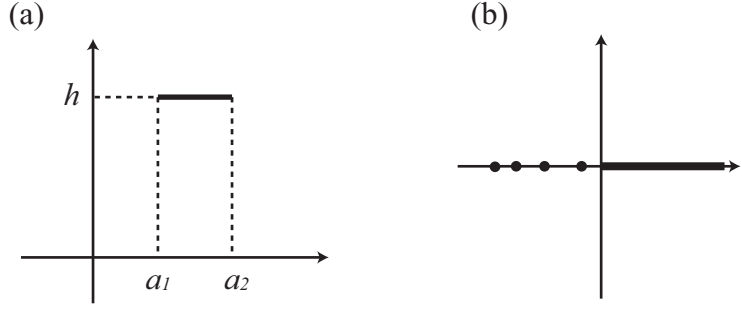}
\caption{(a) A potential $V(x)$ given by Eq.(\ref{4-17}). When $h < 0$, it is well known as a potential well in quantum mechanics.
(b) The spectrum of $T = -\Delta + V$ for $h<0$ in $L^2(\mathbf{R})$-sense, which consists of the continuous
spectrum on the positive real axis and a finite number of eigenvalues on the negative real axis.\label{fig3}}
\end{center}
\end{figure}

For $V$ given by Eq.(\ref{4-17}), the equation $(\lambda -T^\times) \mu = 0$ is of the form
\begin{equation}
0 = \mu'' + (\lambda -V(x)) \mu = \left\{ \begin{array}{ll}
 \mu'' + (\lambda -h) \mu & (a_1 < x < a_2), \\
 \mu'' + \lambda \mu & (\mathrm{otherwise}). \\
\end{array} \right.
\end{equation}
Suppose that $\lambda \neq 0, h$.
A general solution is given by
\begin{equation}
\mu_0 (x) := \left\{ \begin{array}{ll}
\mu_1 (x) = C_1 e^{\sqrt{h-\lambda }x} + C_2 e^{-\sqrt{h-\lambda }x} & (a_1 < x < a_2), \\
\mu_2 (x) = C_3 e^{\sqrt{-\lambda }x} + C_4 e^{-\sqrt{-\lambda }x} & (x < a_1), \\
\mu_3 (x) = C_5 e^{\sqrt{-\lambda }x} + C_6 e^{-\sqrt{-\lambda }x} & (a_2 < x). 
\end{array} \right.
\label{4-19}
\end{equation}
The domain of the dual operator $\Delta^\times$ of the Laplace operator is $\bigcup_{a\geq 1}H^2(\mathbf{R},e^{-2a|x|}dx)$.
Thus we seek a solution in $\bigcup_{a\geq 1} H^2(\mathbf{R}, e^{-2a|x|}dx)$.
We can verify that $\mu_0 (x) \in \bigcup_{a\geq 1} H^2(\mathbf{R}, e^{-2a|x|}dx)$ if and only if
\begin{equation}
\mu_1(a_1) = \mu_2(a_1),\,\, \mu'_1(a_1) = \mu'_2(a_1),\,\, \mu_1(a_2) = \mu_3(a_2),\,\, \mu'_1(a_2) = \mu'_3(a_2). 
\end{equation}
This yields
\begin{eqnarray}
& & \left(
\begin{array}{@{\,}cc@{\,}}
e^{\sqrt{-\lambda }a_1} & e^{-\sqrt{-\lambda }a_1} \\
\sqrt{-\lambda }e^{\sqrt{-\lambda }a_1} & -\sqrt{-\lambda }e^{-\sqrt{-\lambda }a_1}
\end{array}
\right) \left(
\begin{array}{@{\,}c@{\,}}
C_3 \\
C_4
\end{array}
\right) = \left(
\begin{array}{@{\,}c@{\,}}
C_1 e^{\sqrt{h-\lambda }a_1} + C_2 e^{-\sqrt{h-\lambda }a_1} \\
\sqrt{h-\lambda }C_1 e^{\sqrt{h-\lambda }a_1} - \sqrt{h-\lambda } C_2 e^{-\sqrt{h-\lambda }a_1}
\end{array}
\right), \quad \quad \quad \label{4-21}\\
& & \left(
\begin{array}{@{\,}cc@{\,}}
e^{\sqrt{-\lambda }a_2} & e^{-\sqrt{-\lambda }a_2} \\
\sqrt{-\lambda }e^{\sqrt{-\lambda }a_2} & -\sqrt{-\lambda }e^{-\sqrt{-\lambda }a_2}
\end{array}
\right) \left(
\begin{array}{@{\,}c@{\,}}
C_5 \\
C_6
\end{array}
\right) = \left(
\begin{array}{@{\,}c@{\,}}
C_1 e^{\sqrt{h-\lambda }a_2} + C_2 e^{-\sqrt{h-\lambda }a_2} \\
\sqrt{h-\lambda }C_1 e^{\sqrt{h-\lambda }a_2} - \sqrt{h-\lambda } C_2 e^{-\sqrt{h-\lambda }a_2}
\end{array}
\right). \label{4-22}
\end{eqnarray}
Once $C_1$ and $C_2$ are given, $C_3, \cdots ,C_6$ are determined through these equations.
This form of $\mu_0 (x)$ with (\ref{4-21}), (\ref{4-22}) gives a necessary condition for $\mu$ to satisfy the 
equation $(id -  A(\lambda )K^\times)\mu = 0$.

The next purpose is to substitute $\mu_0 (x)$ into Eq.(\ref{3-27}) to determine $C_1, C_2$ and $\lambda $.
Then, we obtain
\begin{eqnarray*}
& & -2\sqrt{-\lambda } (C_1 e^{\sqrt{h-\lambda }x} + C_2 e^{-\sqrt{h-\lambda }x}) \\
&=& h \left( \int^{a_2}_{x}\! e^{-\sqrt{-\lambda } (y-x)}(C_1 e^{\sqrt{h-\lambda }y} + C_2 e^{-\sqrt{h-\lambda }y})dy
 + \int^{x}_{a_1}\! e^{\sqrt{-\lambda } (y-x)}(C_1 e^{\sqrt{h-\lambda }y} + C_2 e^{-\sqrt{h-\lambda }y})dy \right),
\end{eqnarray*}
for $a_1 < x < a_2$.
Comparing the coefficients of $e^{\pm \sqrt{-\lambda }x}$ and $e^{\pm \sqrt{h-\lambda }x}$ in both sides,
we obtain
\begin{eqnarray}
e^{\sqrt{-\lambda }x} : & &  0 = \frac{e^{(\sqrt{h-\lambda }-\sqrt{-\lambda })a_2}}{\sqrt{h-\lambda }-\sqrt{-\lambda }}C_1
 - \frac{e^{-(\sqrt{h-\lambda }+\sqrt{-\lambda })a_2}}{\sqrt{h-\lambda }+\sqrt{-\lambda }}C_2, \label{4-24a}\\
e^{-\sqrt{-\lambda }x} : & & 0 =-\frac{e^{(\sqrt{h-\lambda }+\sqrt{-\lambda })a_1}}{\sqrt{h-\lambda }+\sqrt{-\lambda }}C_1
 + \frac{e^{-(\sqrt{h-\lambda }-\sqrt{-\lambda })a_1}}{\sqrt{h-\lambda }-\sqrt{-\lambda }}C_2,\label{4-24b}
\end{eqnarray}
and
\begin{eqnarray*}
e^{\sqrt{h-\lambda }x} : & & -2\sqrt{-\lambda }C_1 = - \frac{h}{\sqrt{h-\lambda } - \sqrt{-\lambda }}C_1
               + \frac{h}{\sqrt{h-\lambda } + \sqrt{-\lambda }}C_1, \\
e^{-\sqrt{h-\lambda }x} : & & -2\sqrt{-\lambda }C_2 = \frac{h}{\sqrt{h-\lambda } + \sqrt{-\lambda }}C_2
               - \frac{h}{\sqrt{h-\lambda } - \sqrt{-\lambda }}C_2.
\end{eqnarray*}
Note that the last two equations are automatically satisfied.
The condition for the first two equations to have a nontrivial solution $(C_1, C_2)$ is 
\begin{equation}
e^{2\sqrt{h-\lambda }(a_2 - a_1)} = \frac{(\sqrt{h-\lambda } - \sqrt{-\lambda })^4}{h^2}, \quad \lambda \neq 0, h.
\label{4-25}
\end{equation}
Therefore, the generalized eigenvalues on the Riemann surface are given as roots of this equation.
Next, when $\lambda =h$, $\mu_1(x)$ in Eq.(\ref{4-19}) is replaced by $\mu_1(x) = C_1x + C_2$.
By a similar calculation as above, it is proved that $\mu_0(x)$ is a generalized eigenfunction if and only if 
$a_2 - a_1 + 2/\sqrt{-h} = 0$.
However, this is impossible because $a_2 - a_1 > 0$.
Hence, $\lambda =h$ is not a generalized eigenvalue.

\begin{figure}
\begin{center}
\includegraphics{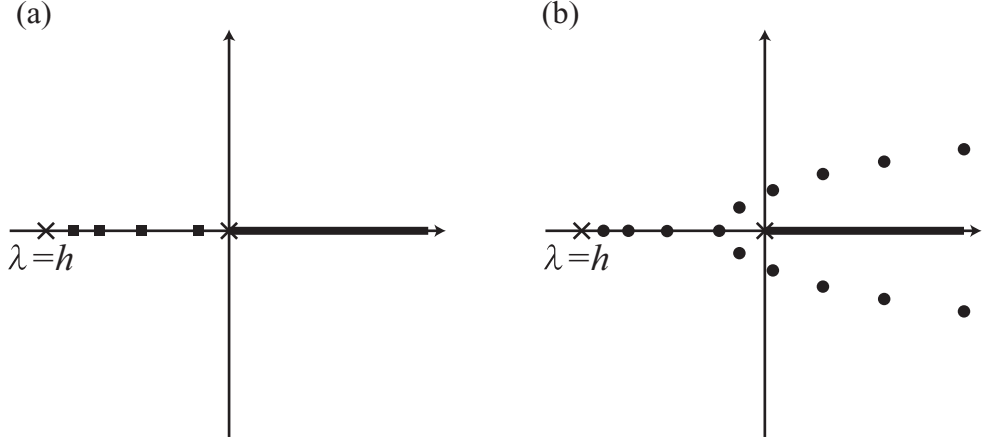}
\caption{Generalized eigenvalues of $T = -\Delta + V$ on (a) the first Riemann sheet and (b) the second Riemann sheet.
The solid lines denote the branch cut. \label{fig4}}
\end{center}
\end{figure}

Note that Eq.(\ref{4-25}) is valid even for $h\in \mathbf{C}$.
Let us investigate roots of Eq.(\ref{4-25}) for $h < 0$ (the case $h>0$ is investigated in a similar manner).
\\[0.2cm]
\textbf{Theorem \thedef.} Suppose that $h<0$.
\\
(I) The generalized spectrum on the first Riemann sheet consists of a finite number of generalized eigenvalues on the 
negative real axis.
\\
(II) The generalized spectrum on the second Riemann sheet consists of
\\
(i) a finite number of generalized eigenvalues on the negative real axis, and
\\
(ii) an infinite number of generalized eigenvalues on both of the upper half plane and the lower half plane. 
They accumulate at infinity along the positive real axis (see Fig.\ref{fig4}).
\\[0.2cm]
This is obtained by an elementary estimate of Eq.(\ref{4-25}) and a proof is omitted.
Put $\lambda = r e^{\sqrt{-1}\theta }$.
We can verify that when $\theta = -\pi$ (the negative real axis on the first Riemann sheet),
generalized eigenvalues are given as roots of one of the equations
\begin{equation}
\tan \Bigl[\frac{1}{2}(a_2-a_1) \sqrt{-r -h}\Bigr] = \frac{\sqrt{r}}{\sqrt{-r-h}}, \quad
\cot \Bigl[\frac{1}{2}(a_2-a_1) \sqrt{-r -h}\Bigr] = -\frac{\sqrt{r }}{\sqrt{-r-h}},
\end{equation}
which have a finite number of roots on the interval $h< \lambda < 0$.
These formulae are well known in quantum mechanics as equations which determine eigenvalues of the Hamiltonian
$-\Delta + V$ (for example, see \cite{Sak}).
Indeed, it is proved in Prop.3.17 of \cite{Chi} that a point $\lambda $ on the first Riemann sheet is an isolated generalized
eigenvalue if and only if it is an isolated eigenvalue of $T$ in the usual sense
(this follows from the fact that $\mathcal{R}_\lambda \circ i = i \circ (\lambda -T)^{-1}$ on the 
 first Riemann sheet).

On the other hand, when $\theta = \pi$ (the negative real axis on the second Riemann sheet), 
Eq.(\ref{4-25}) is reduced to one of the equations
\begin{equation}
\tan \Bigl[\frac{1}{2}(a_2-a_1) \sqrt{-r -h}\Bigr] = -\frac{\sqrt{r }}{\sqrt{-r -h}}, \quad
\cot \Bigl[\frac{1}{2}(a_2-a_1) \sqrt{-r -h}\Bigr] = \frac{\sqrt{r}}{\sqrt{-r -h}},
\end{equation}
which have a finite number of roots on the interval $h< \lambda < 0$.
The theorem (II)-(ii) can be proved by a standard perturbation method.

A straightforward calculation shows that the function $\mathbb{D}(\lambda , \nu_1, \nu_2)$ is given by
\begin{equation}
\mathbb{D}(\lambda , \nu_1, \nu_2)
 = \frac{he^{(\sqrt{h-\lambda } - \sqrt{-\lambda })a_2 - (\sqrt{h-\lambda } - \sqrt{-\lambda })a_1}}
{4\sqrt{-\lambda } (\sqrt{h-\lambda } - \sqrt{-\lambda })^2}
\cdot \left( 1 - \frac{(\sqrt{h-\lambda } - \sqrt{-\lambda })^4}{h^2} e^{-2 \sqrt{h-\lambda } (a_2 - a_1)}\right).
\end{equation}
Thus $\mathbb{D}(\lambda , \nu_1, \nu_2) = 0$ gives the same formula as (\ref{4-25}).
\\

Next, let us consider the one dimensional potential
\begin{equation}
V(x) = - \frac{V_0}{(\cosh x)^2}, \quad V_0\in \mathbf{C},
\end{equation}
which is called the hyperbolic P\"{o}schl-Teller potential \cite{Der}.
This potential satisfies Eq.(\ref{4-11}) for $a=1$ and $V\in G(-\alpha ,\alpha )$ for any $0<\alpha < \pi/2$.
Thus, the generalized resolvent $\mathcal{R}_\lambda i\phi$ of $H+K=-\Delta +V$ is meromorphic on 
$P(1) = \{ \lambda \, | \, -1 <\mathrm{Im}(\sqrt{\lambda }) < 1\}$ when $\phi \in L^2(\mathbf{R}, e^{2|x|}dx)$,
and meromorphic on 
\begin{equation}
\{ \lambda \, | \, -2\pi-2\alpha < \mathrm{arg} (\lambda ) < 2\alpha \},
\end{equation}
when $\phi \in G(-\alpha ,\alpha )$.
In particular, $\mathcal{R}_\lambda i\phi$ is meromorphic on 
\begin{equation}
\hat{\Omega } := \{ \lambda \, | \, -3\pi < \mathrm{arg} (\lambda ) < \pi \},
\end{equation}
if $\phi \in \bigcap _{0<\alpha <\pi/2}G(-\alpha ,\alpha )$.
Let us calculate generalized eigenvalues on the region $\hat{\Omega }$.

Recall that $\mathcal{R}_\lambda $ coincides with $(\lambda -T)^{-1}$ when $-2\pi < \mathrm{arg} (\lambda ) < 0$
(the first Riemann sheet), so that $\lambda $ is a generalized eigenvalue if and only if it is an eigenvalue of 
$-\Delta +V$ in the usual sense.
It is known that eigenvalues $\lambda $ are given by
\begin{equation}
\sqrt{-\lambda } = \pm \sqrt{V_0 + \frac{1}{4}} - n - \frac{1}{2},
\quad \mathrm{Re}(\sqrt{-\lambda }) >0,\,\, n=1,2,\cdots .
\label{3-36}
\end{equation}
It turns out that when $V_0 > 2$, there exists a finite number of eigenvalues on the negative real axis,
and when $V_0<2$, there are no eigenvalues.

\begin{figure}
\begin{center}
\includegraphics{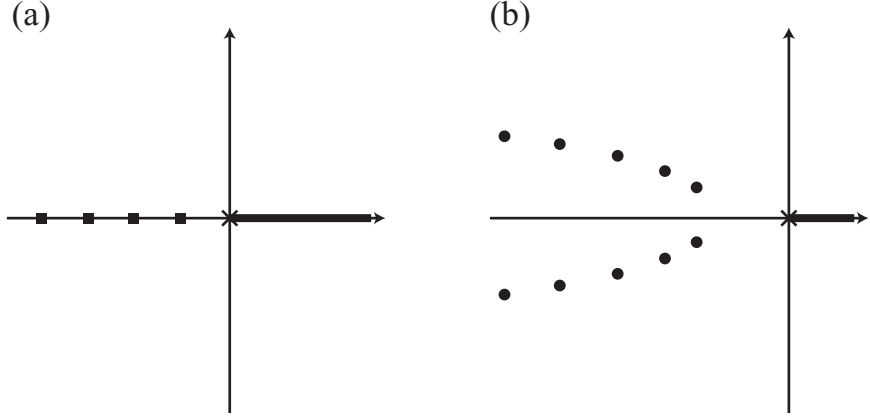}
\caption{Generalized eigenvalues of $T = -\Delta + V$ on (a) the first Riemann sheet for the case $V_0 > 2$,
(b) the second Riemann sheet for the case $V_0 <-1/4$. 
The solid lines denote the branch cut. \label{fig5}}
\end{center}
\end{figure}

In what follows, we suppose that $-3\pi <\mathrm{arg} (\lambda ) < -2\pi$ or 
$0<\mathrm{arg} (\lambda ) < \pi$ to obtain generalized eigenvalues on the second Riemann sheet.
Choose the branch of $\sqrt{-\lambda }$ so that $\mathrm{Re}(\sqrt{-\lambda }) <0$.
\\[0.2cm]
\textbf{Theorem \thedef.}
Generalized eigenvalues on the second Riemann sheet of $\hat{\Omega }$ are given by
\begin{equation}
\sqrt{-\lambda } = \pm \sqrt{V_0 + \frac{1}{4}} - n - \frac{1}{2},
\quad \mathrm{Re}(\sqrt{-\lambda }) <0, \,\, n=1,2,\cdots .
\label{3-37}
\end{equation}
When $V_0 \geq -1/4$, 
there are no generalized eigenvalues (although Eq.(\ref{3-37}) has roots on the ray $\mathrm{arg} (\lambda ) = \pi, -3\pi$).
When $V_0 < -1/4$, there exists infinitely many generalized eigenvalues on the upper
half plane and the lower half plane (see Fig.\ref{fig5}).
Note that Eqs.(\ref{3-36}) and (\ref{3-37}) are valid even when $V_0 \in \mathbf{C}$.
\\

\noindent \textbf{Proof.}
A generalized eigenfunction is defined as a solution of the equation $\mu = A(\lambda )K^\times \mu$,
where $\mu = \mu(x)$ satisfies the differential equation
\begin{equation}
\frac{d^2\mu}{dx^2} + (\lambda + \frac{V_0}{(\cosh x)^2})\mu = 0.
\label{3-38}
\end{equation}
Change the variables by $\rho = (e^{2x} + 1)^{-1}$ and $\mu = (\rho - \rho ^2)^{\sqrt{-\lambda }/2}g$.
Then, the above equation is reduced to the hypergeometric differential equation of the form
\begin{equation}
\rho (1- \rho)g'' + (c-(a+b-1))g' - ab g = 0,
\end{equation}
where
\begin{equation}
a = \sqrt{-\lambda } + \frac{1}{2} - \sqrt{V_0 + \frac{1}{4}},
\quad b = \sqrt{-\lambda } + \frac{1}{2} + \sqrt{V_0 + \frac{1}{4}},
\quad c = \sqrt{-\lambda }+1.
\end{equation}
Hence, a general solution of Eq.(\ref{3-38}) is given by
\begin{eqnarray}
f(\rho):= \mu(x)& =&  A\rho^{\sqrt{-\lambda }/2}(1-\rho)^{\sqrt{-\lambda }/2} F(a,b,c,\rho)\nonumber \\
& &+ B \rho^{-\sqrt{-\lambda }/2}(1-\rho)^{\sqrt{-\lambda }/2} F(a-c+1,b-c+1,2-c,\rho),
\label{3-43}
\end{eqnarray}
where $F(a,b,c, \rho)$ is the hypergeometric function and $A, B\in \mathbf{C}$ are constants.
Note that when $x\to \infty$ and $x\to -\infty$, then $\rho \to 0$ and $\rho \to 1$, respectively.

Fix two positive numbers $\beta < \alpha < \pi/2$.
Put $\rho_{-\alpha } = (e^{2xe^{-\sqrt{-1}\alpha }} + 1)^{-1}$ and define the curve $C$ on $\mathbf{C}$,
which connects $1$ with $0$, to be
\begin{equation}
C = \{ \rho_{-\alpha } = (e^{2xe^{-\sqrt{-1}\alpha }} + 1)^{-1} \, | \, -\infty<x<\infty\}.
\end{equation}
\textbf{Lemma \thedef.} If $\mu$ is a generalized eigenfunction, $f(\rho_{-\alpha })$ is bounded 
as $\rho_{-\alpha } \to 0, 1$ along the curve $C$.
\\[0.2cm]
\textbf{Proof.}
Since $\mu = A(\lambda )K^\times \mu$, $\mu$ is included in the range of $A(\lambda )$.
Then, the estimate (\ref{3-34}) (for $V=1$ and $\theta =-\alpha $) 
proves that $\mu (xe^{-\sqrt{-1}\alpha })$ is bounded as $x\to \pm \infty$.
By the definition of $f$, $f(\rho_{-\alpha }) = \mu (xe^{-\sqrt{-1}\alpha })$. \hfill $\blacksquare$
\\[0.2cm]
\textbf{Lemma \thedef.} 
Suppose $0\leq \mathrm{arg} (\lambda ) < 2\alpha $.
Then, $\rho_{-\alpha }^{\sqrt{-\lambda }} \to 0$ as $\rho_{-\alpha } \to 0$ along the curve $C$.
\\[0.2cm]
\textbf{Proof.}
When $x$ is sufficiently large, we obtain
\begin{equation}
\log \rho_{-\alpha } = -\log (e^{2x e^{-\sqrt{-1}\alpha }} + 1) \sim -\log e^{2x e^{-\sqrt{-1}\alpha }} = -2x e^{-\sqrt{-1}\alpha }.
\end{equation}
Put $\lambda = re^{\sqrt{-1}\mathrm{arg} (\lambda )}$ and 
$\sqrt{-\lambda } = \sqrt{r} e^{\sqrt{-1}(\mathrm{arg}(\lambda ) + \pi)/2}$.
Then,
\begin{eqnarray*}
\rho_{-\alpha }^{\sqrt{-\lambda }} = e^{\sqrt{-\lambda } \log \rho_{-\alpha }} \sim
\exp [-2x\sqrt{r} e^{\sqrt{-1}(\mathrm{arg}(\lambda ) + \pi - 2\alpha )/2}].
\end{eqnarray*}
When $0\leq \mathrm{arg} (\lambda ) < 2\alpha $, this tends to zero as $x\to \infty \, (\rho_{-\alpha } \to 0)$.
 \hfill $\blacksquare$
\\[0.2cm]
\textbf{Lemma \thedef.} When $0\leq \mathrm{arg} (\lambda ) < 2\alpha $, then $B = 0$ and $\Gamma (a) \Gamma (b) = \infty$,
where $\Gamma$ is the gamma function.
\\[0.2cm]
\textbf{Proof.} Due to the definition of the hypergeometric function, we have $F(a,b,c, 0) = 1$.
Then, Eq.(\ref{3-43}), Lemma 3.8 and Lemma 3.9 prove
\begin{equation}
\lim_{\rho_{-\alpha } \to 0} f(\rho_{-\alpha }) = \lim_{\rho_{-\alpha }\to 0} B \rho_{-\alpha }^{-\sqrt{-\lambda }/2} < \infty.
\end{equation}
This yields $B=0$.
Next, we obtain
\begin{equation}
\lim_{\rho_{-\alpha } \to 1} f(\rho_{-\alpha }) 
= \lim_{\rho_{-\alpha }\to 1} A (1-\rho_{-\alpha })^{\sqrt{-\lambda }/2}F(a,b,c, \rho_{-\alpha }).
\end{equation}
It is known that $F(a,b,c, \rho)$ satisfies
\begin{eqnarray}
F(a,b,c, \rho) &=& \frac{\Gamma (c)\Gamma(c-a-b)}{\Gamma(c-a) \Gamma(c-b)} F(a,b, a+b-c+1,1- \rho) \nonumber \\
& +&  (1-\rho)^{c-a-b} \frac{\Gamma (c) \Gamma (a+b-c)}{\Gamma (a) \Gamma (b)} F(c-a, c-b, c-a-b+1, 1-\rho) \quad\quad
\label{3-48}
\end{eqnarray}
for any $\rho$.
This formula is applied to yields
\begin{eqnarray*}
\lim_{\rho_{-\alpha } \to 1} f(\rho_{-\alpha }) =
\lim_{\rho_{-\alpha }\to 1} A (1-\rho_{-\alpha })^{-\sqrt{-\lambda }/2}
\frac{\Gamma (c)\Gamma(a+b-c)}{\Gamma (a) \Gamma (b)} F(c-a, c-b, c-a-b+1, 1-\rho _{-\alpha }).
\end{eqnarray*}
Since $\lim_{\rho_{-\alpha } \to 1} f(\rho_{-\alpha })$ is bounded by Lemma 3.8,
we need $\Gamma (a) \Gamma (b) = \infty$. \hfill $\blacksquare$

The gamma function $\Gamma (x)$ diverges if and only if $x$ is a negative integer.
This proves that Eq.(\ref{3-37}) is a necessary condition for $\lambda $ to be a generalized eigenvalue
when $0\leq \mathrm{arg} (\lambda ) < 2\alpha $.
Since $\alpha  < \pi /2$ is an arbitrary number, the same is true for $0\leq \mathrm{arg} (\lambda ) < \pi$.
The other case $-3\pi< \mathrm{arg} (\lambda ) \leq 2\pi$ is proved in the same way by
replacing $\alpha $ with $-\alpha $.

Finally, let us show that Eq.(\ref{3-37}) is a sufficient condition.
In our situation, the equation $\mu = A(\lambda )K^\times \mu$ is written as
\begin{eqnarray*}
-2\sqrt{-\lambda } \mu (xe^{\sqrt{-1}\theta }) 
&=& e^{\sqrt{-1}\theta } \int^{\infty}_{x}\! e^{-\sqrt{-\lambda }e^{\sqrt{-1}\theta }(y-x)}
V(ye^{\sqrt{-1}\theta })\mu (ye^{\sqrt{-1}\theta })dy \\
& &  + e^{\sqrt{-1}\theta }\int^{x}_{-\infty}\! e^{\sqrt{-\lambda } e^{\sqrt{-1}\theta }(y-x)}
V(ye^{\sqrt{-1}\theta })\mu (ye^{\sqrt{-1}\theta })dy.
\end{eqnarray*}
Note that this is also obtained by an analytic continuation of Eq.(\ref{3-27}).
Putting $\eta _\theta = (e^{2y e^{\sqrt{-1}\theta }} + 1)^{-1}$ provides
\begin{eqnarray*}
-\frac{\sqrt{-\lambda }}{V_0} f(\rho_\theta )
&=& \rho_\theta^{-\sqrt{-\lambda }/2} (1- \rho_\theta)^{\sqrt{-\lambda } /2}
\int^{0}_{\rho_\theta}\! \eta_\theta^{\sqrt{-\lambda }/2} (1- \eta_\theta)^{-\sqrt{-\lambda } /2}
f (\eta_\theta)d\eta_\theta  \\
& &  + \rho_\theta^{\sqrt{-\lambda }/2} (1- \rho_\theta)^{-\sqrt{-\lambda } /2}
\int^{\rho_\theta }_{1}\! \eta_\theta^{-\sqrt{-\lambda }/2} (1- \eta_\theta)^{\sqrt{-\lambda }/2}
f (\eta_\theta)d\eta_\theta ,
\end{eqnarray*}
where the integrals in the right hand side are done along the path $C$.
Since $B= 0$, we obtain
\begin{eqnarray}
-\frac{\sqrt{-\lambda }}{V_0} F(a,b,c, \rho_\theta )
&=& \rho_\theta^{-\sqrt{-\lambda }}
\int^{0}_{\rho_\theta}\! \eta_\theta^{\sqrt{-\lambda }}F (a,b,c, \eta_\theta )d\eta_\theta  \nonumber \\
& &  + (1- \rho_\theta)^{-\sqrt{-\lambda } }
\int^{\rho_\theta }_{1}\! (1- \eta_\theta)^{\sqrt{-\lambda }}F (a,b,c, \eta_\theta )d\eta_\theta .
\label{3-49}
\end{eqnarray}
It is known that the hypergeometric function satisfies
\begin{eqnarray}
& & F(a,b,c, \rho) = (1-\rho)^{c-a-b} F(c-a, c-b, c, \rho), \\
& & F(a,b,c,\rho) = \frac{c-1}{(a-1)(b-1)}\frac{dF}{d\rho}(a-1, b-1, c-1, \rho).
\end{eqnarray}
By using them, Eq.(\ref{3-48}) and $\Gamma (a)\Gamma (b) = 0$, we can verify that Eq.(\ref{3-49}) is reduced to
\begin{eqnarray}
-\frac{\sqrt{-\lambda }}{V_0} F(a,b,c, \rho_\theta )
 = -\frac{\rho_\theta }{c} F(a,b,c+1, \rho_\theta ) - \frac{\sqrt{-\lambda }}{V_0} (1- \rho_\theta)F(a,b,c-1,\rho_\theta).
\end{eqnarray}
Then, the definition $F(a,b,c, \rho) = \sum^\infty_{n=0}(a)_n(b)_n/((c)_n n!) \rho^n$ of the 
hypergeometric function proves that the right hand side of the above equation coincides with the left hand side
for any $\rho_\theta $.
This implies that the function $f(\rho) = \mu (x)$ given by Eq.(\ref{3-43}) satisfies the equation
$\mu = A(\lambda )K^\times \mu$ when $B = 0$ and $\Gamma (a)\Gamma (b) = 0$.
Therefore,  Eq.(\ref{3-37}) is a sufficient condition for $\lambda $ to be a generalized eigenvalue. \hfill $\blacksquare$
\\[0.2cm]
\textbf{Remark.}
For the Schr\"{o}dinger operator $-h^2 \Delta + V$ with the Planck constant $h>0$,
generalized eigenvalues on the second Riemann sheet are given by
\begin{equation}
\sqrt{-\lambda } = \pm \sqrt{V_0 + \frac{h^2}{4}} - h(n + \frac{1}{2}),
\quad \mathrm{Re}(\sqrt{-\lambda }) <0, \,\, n=1,2,\cdots .
\end{equation}
In particular, if $V_0 < 0$ is independent of $h$, we obtain
\begin{equation}
\lambda = -V_0 \pm 2\sqrt{-1}\sqrt{-V_0}(n + \frac{1}{2})h + O(h^2), \quad h\to 0.
\end{equation}
This result coincides with the asymptotic distribution of resonances obtained by \cite{Bri, Sjo}.


\vspace*{0.5cm}
\textbf{Acknowledgements.}

This work was supported by Grant-in-Aid for Young Scientists (B), No.22740069 from MEXT Japan.


\end{document}